\documentclass[11pt]{amsart} 

\usepackage[dvips]{graphics} 
\usepackage{color,amssymb,amsbsy,amsmath,amsfonts,amssymb,
amscd,epsfig,times,graphics}
\newtheorem{theo}{Theorem}[section] 
\newtheorem{prop}[theo]{Proposition}

\newtheorem{lem}[theo]{Lemma} 
\newtheorem{cor}[theo]{Corollary}

\theoremstyle{definition}
\newtheorem{defi}[theo]{Definition}

\numberwithin{equation}{section}

\newcommand{\C}{\mathbb C} 
\newcommand{\R}{\mathbb R} 
 
\newcommand{\B}{\mathbb B}

\newcommand{\partx}{\partial/\partial x}
\begin{document} 
\begin{abstract}
Let $D$ be a $J$-pseudoconvex region in a smooth almost complex
manifold  $(M,J)$ of real dimension four. We construct a local peak $J$-plurisubharmonic function
at every point $p \in bD$ of finite D'Angelo type. 
As applications we give local estimates of the Kobayashi pseudometric, implying the local
Kobayashi hyperbolicity of $D$ at $p$. In case the point $p$ is of D'Angelo type
less than or equal to four, or the approach is nontangential, we
provide sharp estimates of the Kobayashi pseudometric.
\end{abstract} 
\title[Pseudoconvex regions of finite D'Angelo type]
{ Pseudoconvex regions of finite D'Angelo type in four dimensional almost complex manifolds} 
\author{Florian Bertrand}
\address{LATP, C.M.I, 39 rue Joliot-Curie 13453 Marseille cedex 13, FRANCE }
\email{bertrand@cmi.univ-mrs.fr}
\subjclass[2000]{Primary 32Q60, 32T25, 32T40, 32Q45, 32Q65}
\keywords{Almost complex structure, peak plurisubharmonic functions, Kobayashi pseudometric, D'Angelo type}
\maketitle 
{\bf \tableofcontents}

\section*{Introduction} 

Analysis on almost complex manifolds recently  became a fondamental tool in 
symplectic geometry with the work of M.Gromov in \cite{gr}. 
The local existence of pseudoholomorphic discs proved by A.Nijenhuis-W.Woolf 
in their paper \cite{ni-wo}, allows to define the Kobayashi pseudometric, which is crucial for local analysis.

In the present paper we  study the behaviour of the Kobayashi pseudometric of a $J$-pseudoconvex region of finite 
D'Angelo type in an almost complex manifold $\left(M,J\right)$ of dimension four. 
Finite D'Angelo type appeared naturally in complex manifolds when considering the boundary behaviour
of the  $\overline{\partial}$ operator (see \cite{d'a1},\cite{d'a2},\cite{ko},\cite{bl-gr}). 
Moreover on complex manifolds of dimension two, the D'Angelo type unifies many type conditions as the finite regular type.  
Finite regular type  was recently characterized intrinsically by J.-F.Barrault-E.Mazzilli \cite{ba-ma} by means of Lie 
brackets, which generalizes in the non integrable case, a result of T.Bloom-I.Graham \cite{bl-gr}. 

\vskip 0,3cm

Our main result is the construction of a local peak $J$-plurisubharmonic function on pseudoconvex regions 
provided by Theorem A (see also Theorem \ref{proppeak}):  

\vskip 0,3cm

\noindent{\bf Theorem A}. 
{\it Let $D=\{\rho<0\}$ be a domain of finite D'Angelo type in an almost complex 
manifold $\left(M,J\right)$ of  dimension four. We suppose that 
$\rho$ is a $\mathcal{C}^{2}$ defining function of $D$, $J$-plurisubharmonic on a neighborhood of $\overline{D}$. 
Let $p \in \partial D$ be a boundary point. Then  there exists a local peak $J$-plurisubharmonic function at $p$.} 

\vskip 0,3cm

Theorem A  allows to localize pseudoholomorphic discs and to obtain lower estimates 
of the Kobayashi pseudometric which provide the local Kobayashi hyperbolicity of $J$-pseudoconvex regions 
of D'Angelo type $2m$ (Proposition \ref{thm} and Proposition \ref{thmext}). 
As an application we prove the $1/2m$-H\"older extension of biholomorphisms up to the boundary (Proposition \ref{corext}).
In order to obtain sharp lower estimates of the Kobayashi pseudometric similar to those given in 
complex manifolds by D.Catlin \cite{ca} (see also \cite{be2}), we consider a natural scaling method. 
However this reveals the fact that for a domain of finite D'Angelo type greater than four, 
the sequence of almost complex structures obtained  by any polynomial scaling process does not converge generically 
to the standard structure; this is presented in the Appendix.  
This may be related to the fact that finite  D'Angelo type is based on purely complex considerations, 
as the boundary behaviour of the Cauchy-Riemann equations. 
Hence we provide sharp lower estimates of the Kobayashi pseudometric
for a region of finite D'Angelo type four (see also Theorem \ref{theo}):  

\vskip 0,3cm

\noindent{\bf Theorem B}. 
{\it Let $D=\{\rho<0\}$ be a  relatively compact domain  of finite D'Angelo type less than or equal to four 
in an almost complex manifold $\left(M,J\right)$ of dimension four, where $\rho$ is a 
$\mathcal{C}^{2}$ defining function of $D$, $J$-plurisubharmonic on a neighborhood of $\overline{D}$. 
Then there is a positive constant $C$ with the following property: 
for every $p \in D$ and every $v \in T_{p}M$ there exists a diffeomophism $\Phi_{p^*}$ in a neighborhood $U$ of $p$, 
such that: 
\begin{equation*}
K_{\left(D,J\right)}\left(p,v\right)\geq C\left(\frac{|\left(d_p\Phi_{p^*}v\right)_1|}{|\rho\left(p\right)|^{\frac{1}{4}}}+\frac{|\left(d_p\Phi_{p^*}v\right)_{2}|}{|\rho\left(p\right)|}\right).
\end{equation*}}
 
\vskip 0,3cm

We point out that the  approach we use, based on some  renormalization principle of pseudoholomorphic discs, 
gives also a different proof of precise lower estimates obtained by 
H.Gaussier-A.Sukhov in \cite{ga-su} for strictly $J$-pseudoconvex domains in arbitrary dimension. 
As an application of Theorem B, we obtain the (local) complete hyperbolicity of 
$J$-pseudoconvex regions of  D'Angelo type less than or equal to  four (Corollary \ref{corcomp})  and 
we give a Wong-Rosay theorem for regions with noncompact automorphisms group  (Corollary \ref{coraut}). 

\vskip 0,3cm

Finally, in order to obtain  precise  estimates  near a point of arbitrary finite  D'Angelo type,  
we are interested in the nontangential behaviour  of  the Kobayashi pseudometric (see also Theorem \ref{theo3}): 
 
\vskip 0,3cm

\noindent{\bf Theorem C}. 
{\it Let $D=\{\rho<0\}$ be a domain  of finite D'Angelo type  
in an almost complex manifold $\left(M,J\right)$ of dimension four, where $\rho$ 
is a $\mathcal{C}^{2}$ defining function of $D$, $J$-plurisubharmonic on a neighborhood of $\overline{D}$. 
Let $q \in \partial D$ be a boundary point of D'Angelo type $2m$ 
and let $\Lambda \subset D$ be a cone with vertex at $q$ and axis the inward normal axis.  
Then there exists a positive constant $C$ such that for every $p \in D\cap \Lambda$ and every $v=v_n+v_t \in T_{p}M$: 
\begin{equation*}
K_{\left(D,J\right)}\left(p,v\right)\geq C\left(\frac{|v_n|}{|\rho\left(p\right)|^{\frac{1}{2m}}}+
\frac{|v_{t}|}{|\rho\left(p\right)|}\right),
\end{equation*}
where $v_n$ and $v_t$ are the normal and the tangential parts of $v$ with respect to $q$.} 
\vskip 0,3cm

\section{Preliminaries}
We denote by $\Delta$ the unit disc of $\C$ and by $\Delta_{r}$ 
the disc of $\C$ centered at the origin of radius $r>0$.
\subsection{Almost complex manifolds and pseudoholomorphic discs}

An {\it almost complex structure} $J$ on a real smooth manifold $M$ is a $\left(1,1\right)$ tensor field
 which  satisfies $J^{2}=-Id$. We suppose that $J$ is smooth.
The pair $\left(M,J\right)$ is called an {\it almost complex manifold}. We denote by $J_{st}$
the standard integrable structure on $\C^{n}$ for every $n$.
A differentiable map $f:\left(M',J'\right) \longrightarrow \left(M,J\right)$ beetwen two almost complex manifolds is said to be 
 {\it $\left(J',J\right)$-holomorphic}  if $J\left(f\left(p\right)\right)\circ d_{p}f=d_{p}f\circ J'\left(p\right),$ 
for every $p \in M'$. In case  $M'=\Delta \subset \C$, such a map is called a {\it pseudoholomorphic disc}.  
If $f:\left(M,J\right)\longrightarrow M'$ is a diffeomorphism, we define an almost complex structure, $f_{*}J$,  
on $M'$ as the
 direct image of $J$ by $f$ :
$$f_{*}J\left(q\right):=d_{f^{-1}\left(q\right)}f\circ J\left(f^{-1}\left(q\right)\right)\circ d_{q}f^{-1},$$ 
 for every  $q \in M'$.

The following lemma (see \cite{ga-su}) states that locally any almost
 complex manifold can be seen as the unit ball of 
$\C^n$ endowed with a small smooth pertubation of the standard
 integrable structure $J_{st}$. 
\begin{lem}\label{ilemloc}
Let $\left(M,J\right)$ be an almost complex manifold, with $J$ of class
 $\mathcal{C}^{k}$, $k\geq 0$. 
Then for every point $p \in M$ and every $\lambda_0 > 0$ there exist a neighborhood $U$ of $p$ and
 a coordinate diffeomorphism $z: U \rightarrow \B$ centered a $p$ (ie $z(p) =0$) such that  the
direct image of $J$ satisfies $z_{*}J\left(0\right) = J_{st}$ and
$||z_*\left(J\right) - J_{st}||_{\mathcal{C}^k\left(\bar {\B}\right)}
 \leq \lambda_0$.
\end{lem}
This is simply done by considering a local chart $z: U \rightarrow \B$ centered a $p$ (ie $z(p) =0$), composing it 
with a linear diffeomorphism to insure $z_{*}J\left(0\right) = J_{st}$ and  dilating coordinates.

So let  $J$ be an almost complex structure defined in  a neighborhood $U$ of the origin in $\R^{2n}$, and such 
that $J$ is sufficiently closed to the standard structure in uniform norm on the closure $\overline{U}$ of $U$. 
The $J$-holomorphy equation for a pseudoholomorphic disc 
$u : \Delta \rightarrow U \subseteq \R^{2n}$ is given by  
\begin{equation}\label{eqholo}
\frac{\partial u}{\partial y}-J\left(u\right)\frac{\partial u}{\partial x}=0.
\end{equation}

According to \cite{ni-wo}, for every $p \in M$, there is a neighborhood $V$ of zero in $T_{p}M$, such that for every $v \in V$, 
there is a $J$-holomorphic disc $u$ satisfying $u\left(0\right)=p$ and $d_{0}u\left(\partx\right)=v$.

\subsection{Levi geometry}

Let $\rho$ be a $\mathcal{C}^2$ real valued function on
a smooth almost complex manifold  $\left(M,J\right).$
We denote by $d^c_J\rho$ the differential 
form defined by 
$$d^c_J\rho\left(v\right):=-d\rho\left(Jv\right),$$
where  $v$ is a section of $TM$. 
The {\it Levi form} of $\rho$ at a point $p\in M$ and a vector 
$v \in T_pM$ is defined by
\begin{equation*}
\mathcal{L}_J\rho\left(p,v\right):=d\left(d^c_J\rho\right)(p)\left(v,J(p)v\right)=dd^c_J\rho(p)\left(v,J(p)v\right).
\end{equation*}
In case $(M,J)=(\C^n,J_{st})$, then $\mathcal{L}_{J_{st}}\rho$ is, up to a positive multiplicative constant, 
the usual standard Levi form :
\begin{equation*}
\displaystyle \mathcal{L}_{J_{st}}\rho(p,v)=4 \sum \frac{\partial^2\rho}{\partial z_j\partial \overline{z_k}}v_j\overline{v_k}.
\end{equation*}

\vspace{0,3cm}

We investigate now how close is the Levi form with respect to $J$ from the standard Levi form. 
For $p \in M$ and $v\in T_pM$, we easily get :
\begin{equation}\label{ieqlevi}
\mathcal{L}_J\rho\left(p,v\right)=\displaystyle \mathcal{L}_{J_{st}}\rho(p,v)+d(d^c_{J}-d^c_{J_{st}})\rho(p)(v,J(p)v)+
dd^c_{J_{st}}\rho(p)(v,J(p)-J_{st})v).
\end{equation}
In local coordinates $(t_1,t_2,\cdots,t_{2n})$ of $\R^{2n}$, (\ref{ieqlevi}) may be written as follows
\begin{eqnarray}\label{ieqlevi2}
\mathcal{L}_J\rho\left(p,v\right)&=&\displaystyle \mathcal{L}_{J_{st}}\rho(p,v)+{}^tv(A-{}^tA)J(p)v+
{ }^t(J(p)-J_{st})vDJ_{st}v+ \nonumber\\
& & { }^t(J(p)-J_{st})vD(J(p)-J_{st})v
\end{eqnarray}
where 
$$A:=\left(\sum_i\frac{\partial u}{\partial t_i}\frac{\partial
    J_{i,j}}{\partial t_k}\right)_{1\le j,k\le 2n}\quad \mathrm{and}\quad
D:=\left(\frac{\partial^2u}{\partial t_j\partial
    t_k}\right)_{1\le j,k\le 2n}.$$

\vspace{0,3cm}

 Let $f$ be a $(J',J)$-biholomorphism 
from $\left(M',J'\right)$ to $\left(M,J\right)$. Then for every $p\in M$ and every $v\in T_pM$:
$$\mathcal{L}_{J'}\rho\left(p,v\right)=\mathcal{L}_{J}\rho\circ
 f^{-1}\left(f\left(p\right),d_pf\left(v\right)\right).$$
This expresses the invariance of the Levi form under pseudobiholomorphisms.

The next proposition is useful in order to compute the Levi form (see \cite{de}, \cite{ha} and \cite{iv-ro}).   
\begin{prop}\label{proplevi}\mbox{ }
Let $p\in M$ and $v\in T_pM$. Then  
$$\mathcal{L}_J\rho\left(p,v\right)=\Delta \left(\rho \circ u\right)
 \left(0\right),$$ 
where $u : \Delta \rightarrow \left(M,J\right)$ is any $J$-holomorphic
 disc satisfying 
$u\left(0\right)=p$ and $d_0u\left(\partial/\partial_x\right)=v$.
\end{prop}

Proposition \ref{proplevi} leads to the following proposition-definition : 
\begin{prop}\label{iproplevi2}
The two statements  are equivalent :
\begin{enumerate}
\item $\rho \circ u$ is subharmonic for any $J$-holomorphic disc $u :
 \Delta \rightarrow M$.
\item $\mathcal{L}_{J}\rho(p,v)\geq 0$ for every $p \in M$ and every $v
 \in T_pM$. 
\end{enumerate}
\end{prop}
If one of the previous statements is satisfied we say that $\rho$ is
 {\it $J$-plurisubharmonic}. We say that $\rho$ is {\it strictly $J$-plurisubharmonic} if
  $\mathcal{L}_{J}\rho(p,v)$ is positive  
for any $p \in M$ and any $v \in T_pM\setminus\{0\}$. 
$J$-plurisubharmonic functions play a very important role in almost complex geometry :  they give attraction and 
localization properties for pseudoholomorphic discs. For this reason the 
construction of $J$-plurisubharmonic functions is crucial. 

\vspace{0,3cm}

Similarly to the integrable case, one may define the notion of pseudoconvexity in almost complex manifolds.
Let $D$ be a domain in $\left(M,J\right)$. We denote by $T^{J}\partial D:=T\partial D\cap JT\partial D$ the
 $J$-invariant subbundle of $T\partial D.$

\begin{defi}\mbox{ }

\begin{enumerate}
\item The domain $D$ is   $J$-pseudoconvex (resp. it strictly
 $J$-pseudoconvex) 
if  $\mathcal{L}_{J}\rho(p,v)\geq 0$ (resp. $>0$) for any 
$p \in \partial D$ and $v \in T^J_p\partial D$ (resp. $v \in
 T^J_p\partial D \setminus\{0\}$).
\item A $J$-pseudoconvex region is a domain $D=\{\rho<0\}$ where $\rho$
 is a $\mathcal{C}^2$ defining function,
$J$-plurisubharmonic on a neighborhood of $\overline{D}$. 
\end{enumerate}
\end{defi}
We recall that a defining function for $D$ satisfies $d\rho\neq 0$ on $\partial D$.

The following Lemma is useful in order to compute the Levi form of some functions.  

\begin{lem}\label{lemlev}
Assume that $J$ is a diagonal  almost complex structure on $\R^4$ that coincides with 
the standard structure $J$ on $\C\times \{0\}$. To fix notations we suppose that its matricial representation is given by :
\begin{equation*}
J=\left(\begin{array}{ccccc} 
a_1& b_1 & 0 & 0 \\
c_1 & -a_1 & 0 & 0\\
0  & 0 & a_2 & b_2\\
0 & 0 & c_2 & -a_2\\ 
\end{array}\right).
\end{equation*}
Then 
the Levi form of some smooth real valued function $f$ at a point $z=(z_1,z_2)$ and $v=\left(1,0,0,0\right)$ is equal to  
\begin{eqnarray*}\label{eqlevi}
\mathcal{L}_{J}f\left(z,v\right)= -c_1\Delta_1 f+O\left(|z_2|\right)).
\end{eqnarray*}
where $\displaystyle \Delta_1 f:=\frac{\partial^2 f}{\partial x_1\partial x_1}+\frac{\partial^2 f}{\partial y_1\partial y_1}$.

\end{lem}
\begin{proof}
Let us  compute the Levi form of some smooth real valued function $f$ at a point $z=(z_1,z_2)$ and $v=\left(1,0,0,0\right)$ :
\begin{eqnarray*}
c_1^{-1}\mathcal{L}_{J}f\left(z,v\right)&=&-\Delta_1 f+\left[
-2\frac{\partial^2 f}{\partial x_1\partial y_1}a_1+\frac{\partial^2 f}{\partial x_1\partial x_1}\left(1+b_1\right)+
\frac{\partial^2 f}{\partial y_1\partial y_1}\left(c_1-1\right)\right]+ \nonumber \\
& & 
\frac{\partial f}{\partial x_1}\left[\frac{\partial b_1}{\partial x_1}-\frac{\partial a_1}{\partial y_1}\right]+ 
\frac{\partial f}{\partial y_1}\left[\frac{\partial a_1}{\partial x_1}+\frac{\partial c_1}{\partial y_1}\right]\nonumber \\
&=&-\Delta_1 f+\left[-2\frac{\partial^2 f}{\partial x_1\partial y_1}O\left(|z_2|\right)+
\frac{\partial^2 f}{\partial x_1\partial x_1}O\left(|z_2|\right)+\frac{\partial^2 f}
{\partial y_1\partial y_1}O\left(|z_2|\right)\right]+\nonumber \\
& & \frac{\partial f}{\partial x_1}O\left(|z_2|\right)+\frac{\partial f}{\partial y_1}O\left(|z_2|\right)\nonumber \\
&= & -\Delta_1 f+O\left(|z_2|\right).
\end{eqnarray*}

\end{proof}

\section{Construction of a local peak plurisubharmonic function}

This section is devoted to the proof of Theorem A (see Theorem \ref{proppeak}).

\subsection{Pseudoconvex regions of finite D'Angelo  type}

In this subsection we  describe  a pseudonconvex region on a neighborhood of a boundary point of finite D'Angelo type.
We point out that all our considerations are purely local. 
Assume that $D=\{\rho<0\}$ is a $J$-pseudoconvex region in $\C^{2}$ and 
that the structure $J$ is defined on a fixed neighborhood $U$ of $\overline{D}$.
We suppose that the origin  is a boundary point of $D$.  
\begin{defi}\label{defco} 
Let $u~: \left(\Delta,0\right)\rightarrow \left(\R^{4},0,J\right) $ be a $J$-holomorphic disc satisfying $u\left(0\right)=0$.
The   order of contact $\delta_0\left(\partial D,u\right)$ with $\partial D$ at the origin is the degree of the first term 
in the Taylor expansion of 
$\rho \circ u$. 
We denote by $\delta\left(u\right)$ the multiplicity of $u$ at the origin. 
\end{defi}

We now define the D'Angelo  type and the regular type of the real hypersurface $\partial D$ at the origin.  

\begin{defi}\label{deftyp}\mbox{ }

 \begin{enumerate}
\item The   D'Angelo type  of $\partial D$ at the origin is defined by:
$$\Delta^1\left(\partial D,0\right):=\sup \left\{\frac{\delta_0\left(\partial D,u\right)}{\delta\left(u\right)}, \mbox{ } 
u:\Delta\rightarrow \left(\R^{4},J\right) \mbox{  $J$-holomorphic },  u\left(0\right)=0 \right\}.$$
The point $0$ is a point of finite D'Angelo type $2m$ if $\Delta^1\left(\partial D,0\right)=2m<+\infty$.
\item The   regular type of $\partial D$ at origin is defined by:
\begin{eqnarray*}
\Delta^1_{{\rm reg}}\left(\partial D,0\right)&:=&\sup \{\delta_0\left(\partial D,u\right), \mbox{ } 
u:\Delta\rightarrow \left(\R^{4},J\right) \mbox{  $J$-holomorphic },\\ 
& & \hspace{6,5cm} u\left(0\right)=0, d_0u\neq 0 \}.
\end{eqnarray*}
\end{enumerate}
\end{defi}
Since the  regular type of $\partial D$ at the origin consists in considering only regular discs we have:
\begin{equation}\label{eqcomp}
\Delta^1_{{\rm reg}}\left(\partial D,0\right) \leq \Delta^1\left(\partial D,0\right).
\end{equation}
The  type condition  as defined in part 1 of Definition \ref{deftyp} was introduced by 
J.-P.D'Angelo \cite{d'a1}, \cite{d'a2} who proved that this coincides  with the regular type in complex manifolds
of dimension two.  
After Proposition \ref{prop1}, we will also prove that the D'Angelo type and  the regular type coincide in 
four dimensional almost complex manifolds (see Proposition \ref{proptyp}).  

\vspace{0,7cm}

We suppose that the origin is a point of finite regular type. Then let $u~: \Delta \rightarrow \R^4$ be a regular 
$J$-holomorphic disc of maximal contact order $2m$. We choose coordinates such that $u$ is given by 
$u\left(\zeta\right)=\left(\zeta,0\right)$, $J\left(z_1,0\right)=J_{st}$ and such that the complex tangent space 
$T_0\partial D \cap J(0)T_0 \partial D$ is equal to $\{z_2=0\}$. 
Then by considering the family of vectors $\left(1,0\right)$ at base points $\left(0,t\right)$ for $t\neq 0$ small enough, we obtain a family of $J$ 
holomorphic discs $u_{t}$ such that $u_t\left(0\right)=\left(0,t\right)$ and $d_0u_t\left(\partial/\partial_x\right)=\left(0,1\right)$. Due to the parameters dependance of the 
solution to the $J$-holomorphy equation (\ref{eqholo}), we straighten these discs into the complex lines $\{z_2=t\}$. 
We then consider a transversal foliation 
by $J$-holomorphic discs and straighten these lines into $\{z_1=c\}$. In these new coordinates still denoted by $z$, 
the matricial representation of $J$ is diagonal: 
\begin{equation}\label{eqstr}
J=\left(\begin{array}{ccccc} 
a_1& b_1 & 0 & 0 \\
c_1 & -a_1 & 0 & 0\\
0  & 0 & a_2 & b_2\\
0 & 0 & c_2 & -a_2\\ 
\end{array}\right).
\end{equation}
Since  $J\left(z_1,0\right)=J_{st}$ we have 
\begin{equation}\label{eqstrO}
J=J_{st}+O\left(|z_2|\right).
\end{equation} 

In the next fundamental proposition we describe   precisely the local  expression of  the  defining function $\rho$.
\begin{prop}\label{prop1}
The $J$-plurisubharmonic defining function for the domain $D$ has the following local expression:
$$\rho=\Re ez_2+H_{2m}\left(z_1,\overline{z_1}\right)+\widetilde{H}(z_1,z_2)+O\left(|z_1|^{2m+1}+|z_2||z_1|^m+|z_2|^2\right)$$ 
where  $H_{2m}$ is a homogeneous polynomial of degree $2m$, subharmonic which  is not harmonic and 
$$\widetilde{H}(z_1,z_2)=\Re e \displaystyle \sum_{k=1}^{m-1} \rho_{k} z_1^k z_2.$$
\end{prop}

\begin{proof}
Since $T_0\partial D \cap J(0)T_0\partial D =\{z_2=0\}$, we have 
$$\rho=\Re ez_2+O(\|z\|^{2}).$$ 
Moreover the disc $\zeta \mapsto \left(\zeta,0\right)$ being a regular $J$-holomorphic disc of maximal contact order $2m$, 
the defining function $\rho$ has the following local expression: 
$$\rho=\Re ez_2+H_{2m}\left(z_1,\overline{z_1}\right)+O\left(|z_1|^{2m+1}+|z_2|\|z\|\right),$$ 
where $H_{2m}$ is a homogeneous polynomial of degree $2m$.
\vspace{0,7cm}

We prove that the polynomial $H_{2m}$ is subharmonic  using a standard dilation argument. 
Consider the non-isotropic dilation  of $\C^2$
$$\Lambda_\delta\left(z_1,z_2\right):=\left(\delta^{-\frac{1}{2m}}z_1,\delta^{-1}z_2\right).$$ 
Due to Proposition \ref{proplevi}, the domain
$$\Lambda_\delta\left(D\right)=\{\delta^{-1}\left(\rho\circ \Lambda_{\delta}^{-1}\left(z_1,z_2\right)\right)<0\}$$ 
is $\left(\Lambda_\delta\right)_*\left(J\right)$-pseudoconvex.  
Moreover $\Lambda_\delta\left(D\right)$ converges in the sense of local Hausdorff set convergence to 
$$\tilde{D}:=\{Re\left(z_2\right) +H_{2m}\left(z_1,\overline{z_1}\right)<0\},$$ 
as $\delta$ tends to zero and the sequence of structures $\left(\Lambda_\delta\right)_*J$ 
converges to the standard structure $J_{st}$. 
It follows that the limit domain  $\tilde{D}$ is $J_{st}$-pseudoconvex implying that $H_{2m}$ is subharmonic.

\vspace{0,7cm}

Now we prove $H_{2m}$ that contains a nonharmonic part. 
By contradiction, we assume that  $H_{2m}$ is  harmonic. Then $H_{2m}$ can be written $\Re e z_1^{2m}$. 
According to Proposition 1.1 of \cite{iv-ro}, and since the structure $J$ is smooth there exists, for a sufficiently small  $\lambda>0$,  
a pseudoholomorphic disc $u~: \Delta \rightarrow (\R^4,J)$  such that:
$$ \left\{  
\begin{array}{lll}  
\displaystyle  u\left(0\right)&=& 0 \\
\\
\displaystyle  \frac{\partial u}{\partial x}\left(0\right) & = & \left(\lambda^{\frac{1}{2m}},0,0,0\right)\\
\\
\displaystyle  \frac{\partial^{k} u}{\partial x^k}\left(0\right) & = & \left(0,0,0,0\right), \mbox{ for  } 1<k<2m \\
\\
\displaystyle  \frac{\partial^{2m} u}{\partial x^{2m}}\left(0\right) & = & \left(0,0,-\lambda\left(2m\right)!,0\right).\\
\end{array} \right.$$
We  prove that the contact order of  such a regular disc $u$ is greater than $2m$ which contradicts the fact that $D$ is of regular type $2m$.
We denote by $[\rho \circ u]_{2m}$ the homogeneous part of degree $2m$ in the Taylor expansion of $\rho \circ u$ at the origin: 
$$[\rho \circ u]_{2m}\left(x,y\right)=\displaystyle \sum_{k=0}^{2m}a_kx^{k}y^{2m-k}.$$
Let us  prove that 
$\displaystyle  a_k=\frac{\partial^{k} }{\partial x^k}\frac{\partial^{2m-k} }{\partial y^{2m-k}}\rho \circ u \left(0\right)$ is 
equal to zero for each $0\leq k \leq 2m$. 

\vspace{0,2cm}

For $a_{2m}$,  we have:
\begin{eqnarray*}
\frac{\partial^{2m} }{\partial x^{2m}}\rho \circ u \left(0\right)& = & \Re e \frac{\partial^{2m} }{\partial x^{2m}}u_2\left(0\right)+
\Re e\frac{\partial^{2m} }{\partial x^{2m}}u_1^{2m}\left(0\right) \\
&&\\
&=&  -\lambda\left(2m\right)!+\Re e\frac{\partial^{2m} }{\partial x^{2m}}u_1^{2m}\left(0\right).           
\end{eqnarray*}
Since $u_1\left(0\right)=0$, it follows that the only non vanishing term in  
$\displaystyle  \Re e\frac{\partial^{2m} }{\partial x^{2m}}u_1^{2m}\left(0\right)$ is 
$$\left(2m\right)!\Re e\left(\frac{\partial u_1 }{\partial x}\left(0\right)\right)^{2m}=\lambda\left(2m\right)!.$$ 
This proves that $a_{2m}=0.$ 

\vspace{0,2cm}

Then let $0 \leq k< 2m$:
\begin{eqnarray*}
\frac{\partial^{k} }{\partial x^k}\frac{\partial^{2m-k} }{\partial y^{2m-k}}\rho \circ u \left(0\right)& = & 
\Re e \frac{\partial^{k} }{\partial x^k}\frac{\partial^{2m-k}}{\partial y^{2m-k}}u_2\left(0\right)
+\Re e\frac{\partial^{k} }{\partial x^k}\frac{\partial^{2m-k}}{\partial y^{2m-k}}u_1^{2m}\left(0\right). 
\end{eqnarray*}
For the same reason as previously, the only term to consider in 
$\displaystyle  \Re e\frac{\partial^{k} }{\partial x^k}\frac{\partial^{2m-k}}{\partial y^{2m-k}}u_1^{2m}\left(0\right)$ is 
$$\left(2m\right)!\Re e\left(\frac{\partial }{\partial x}u_1\left(0\right)\right)^k
\left(\frac{\partial}{\partial y}u_1\left(0\right)\right)^{2m-k}=
\lambda^{\frac{k}{2m}}\left(2m\right)!\Re e\left(\frac{\partial u_1}{\partial y}\left(0\right)\right)^{2m-k}.$$
Then, since $u$ is $J$-holomorphic, it satisfies  the diagonal $J$-holomorphy equation: 
$$\frac{\partial u_l}{\partial y}=J_l\left(u\right)\frac{\partial u_l}{\partial x}, $$
for $l=1,2$, where 
$$J_l=\left(\begin{array}{cc} 
a_l& b_l \\
c_l & -a_l \\
\end{array}\right)  \mbox{ }\mbox{ (see (\ref{eqstr}) for notations).}$$
It follows that 
\begin{eqnarray*}
\lambda^{\frac{k}{2m}}\left(2m\right)!\Re e\left(\frac{\partial u_1}{\partial y}\left(0\right)\right)^{2m-k} &= &
\lambda^{\frac{k}{2m}}\left(2m\right)!\Re e\left(J_1\left(u\left(0\right)\right)\frac{\partial u_1}{\partial x}
\left(0\right)\right)^{2m-k}\\
&&\\
&=& \lambda\left(2m\right)!\Re e\left(i\right)^{2m-k}.\\
\end{eqnarray*}
Moreover due to the condition $\displaystyle  \frac{\partial^{k} u_2}{\partial x^k}\left(0\right)=\left(0,0\right)$,  for $ 1\leq k<2m$, it follows that 
 the only part  we need to consider in $\displaystyle  \frac{\partial^{2m-k}}{\partial y^{2m-k}}u_2\left(0\right)$ is 
$\displaystyle  J_2\left(u\right)\frac{\partial }{\partial x}
\frac{\partial^{2m-k-1}}{\partial y^{2m-k-1}}u_2\left(0\right)$ and by induction 
$\displaystyle  \left(J_2\left(u\right)\right)^{2m-k}\frac{\partial^{2m-k}}{\partial x^{2m-k}}u_2\left(0\right)$. 
Finally 
\begin{eqnarray*}
\Re e \frac{\partial^{k} }{\partial x^k}\frac{\partial^{2m-k}}{\partial y^{2m-k}}u_2\left(0\right) & = & 
\Re e\left(J_2\left(u\left(0\right)\right)\right)^{2m-k}\frac{\partial^{2m}u_{2}}{\partial x^{2m}}\left(0\right)\\
&&\\
& = & -\lambda\left(2m\right)!\Re e\left(J_2\left(u\left(0\right)\right)^{2m-k}\left(1,0\right)\right)\\
&&\\
&=& -\lambda\left(2m\right)!\Re e\left(i\right)^{2m-k}.\\
\end{eqnarray*}
This proves that the homogeneous part   $[\rho \circ u]_{2m}$ is equal to zero. 

For smaller order terms it is a direct consequence of $u\left(0\right)=0$ and 
$\displaystyle  \frac{\partial^{k} u}{\partial x^k}\left(0\right)  =  \left(0,0,0,0\right),$ for $1<k<2m$.

\vspace{0,7cm}

It remains to prove there are no term  $\Re e \rho_k z_1^k\overline{z_2}$ with $k<m$ in the defining function $\rho$. 
This is done by contradiction and by computing the Levi form of $\rho$ at a point 
$z_0=\left(z_1,0\right)$ and at a vector  $v=\left(X_1,0,X_2,0\right)$. Assume that 
$$\rho=\Re ez_2+H_{2m}\left(z_1,\overline{z_1}\right)+\widetilde{H}(z_1,z_2)+ \Re e \rho_{\overline{k}} z_1^k\overline{z_2} +O\left(|z_1|^{2m+1}+|z_2||z_1|^{k+1}+|z_2|^2\right),$$
with $k<m$. Replacing $z_1$ by $\displaystyle (\rho_{\overline{k}})^{\frac{1}{k}} z_1$ if necessary, we suppose $\rho_{\overline{k}}=1$.
  
The Levi form of $\Re e z_2$ at a point $z_0=\left(z_1,0\right)$ and at a vector  $v=\left(X_1,0,X_2,0\right)$ is equal to 
\begin{eqnarray*}
\mathcal{L}_{J}\Re e z_2\left(z_0,v\right)&=& \left[\left(a_1-a_2\right)(z_0)\frac{\partial a_2}{\partial x_1}(z_0)+
c_1(z_0)\frac{\partial a_2}{\partial y_1}(z_0)
-c_2(z_0)\frac{\partial b_2}{\partial x_1}(z_0)\right]X_1X_2+\\
&&\\
& & c_2(z_0)\left[\frac{\partial a_2}{\partial y_2}(z_0)-\frac{\partial b_2}{\partial x_2}(z_0)\right]X_2^2.
\end{eqnarray*}
Due to (\ref{eqstrO}) we have 
\begin{equation*}
\left\{
\begin{array}{lllll}  
a_1\left(z_0\right)&=&a_2\left(z_0\right)&=&0,\\
\\
c_2\left(z_0\right)&=&1,&&\\
\\
\displaystyle \frac{\partial a_2}{\partial y_1}\left(z_0\right)&=&\displaystyle
\frac{\partial b_2}{\partial x_1}\left(z_0\right)&=&0.
\end{array}
\right.
\end{equation*}
%$$a_1\left(z_0\right)=a_2\left(z_0\right)=0,  \mbox{ }  c_2\left(z_0\right)=1,  \mbox{ } 
% \frac{\partial a_2}{\partial y_1}\left(z_0\right)=0, \mbox{ } \mbox{ and } \mbox{ }  
%\frac{\partial b_2}{\partial x_1}\left(z_0\right)= 0.$$
So the Levi form of $\Re e z_2$ at  $z_0=\left(x_1,0,0,0\right)$ and at a vector  $v=\left(X_1,0,X_2,0\right)$ is 
$$\mathcal{L}_{J}\Re e z_2\left(z_0,v\right)= \left[\frac{\partial a_2}{\partial y_2}(z_0)-\frac{\partial b_2}{\partial x_2}(z_0)\right]X_2^2.$$

According to Lemma \ref{lemlev}, the Levi form of $H_{2m}+O(|z_1|^{2m+1})$ at $z_0$ and $v_1=(X_1,0,X_2,0)$ is equal to  
\begin{eqnarray*}
\mathcal{L}_{J}(H_{2m}+O(|z_1|^{2m+1})) \left(z_0,v\right)& =& \Delta\left(H_{2m}+O(|z_1|^{2m+1})\right) X_1^2+O(|z_1|^{2m-1})X_1X_2.
\end{eqnarray*}

According to the fact  that the Levi form for the standard structure of $\widetilde{H}(z_1,z_2)$ is identically equal to zero, 
and due to (\ref{ieqlevi2}) and to (\ref{eqstrO}), it follows that the Levi form of $\widetilde{H}(z_1,z_2)$ at $z_0$ is equal to 
\begin{eqnarray*}
\mathcal{L}_{J}\widetilde{H} \left(z_0,v\right)& =& O(|z_1|)X_2^2.
\end{eqnarray*}

Now the Levi form of $O(|z_2|^2)$ is equal to 
\begin{eqnarray*}
\mathcal{L}_{J}O(|z_2|^2)\left(z_0,v\right)& =& O(1)X_2^2.
\end{eqnarray*}

And the Levi form of $\Re e z_1^k\overline{z_2}$ is equal 
\begin{eqnarray*}
\mathcal{L}_{J}\Re e z_1^k\overline{z_2}\left(z_0,v\right)& =& (k\Re e z_1^{k-1})X_1X_2+O(|z_1|^k)X_2^2.
\end{eqnarray*}

Finally the Levi form of the defining function $\rho$ at a point $z_0=\left(z1,0\right)$ and at a vector  
$v=\left(X_1,0,X_2,0\right)$ is equal to: 
\begin{eqnarray*}
\mathcal{L}_{J}\rho\left(z_0,v\right)& =& O\left(|z_1|^{2m-2}\right)X_1^2+\left[4k\Re e z_1^{k-1}+O(|z_1|^{2m-1})\right]X_1X_2 \\
&&\\
& & +\left[\frac{\partial a_2}{\partial y_2}(z_0)-\frac{\partial b_2}{\partial x_2}(0)+
O(1)+O\left(|z_1|\right)\right]X_2^2.\\
\end{eqnarray*}
It follows that since $k<m$ there are $z_1$, $X_1$ and $X_2$ such that $\mathcal{L}_{J}\rho\left(z_0,v\right)$ is negative, providing a contradiction.

\end{proof}

Now we prove that the D'Angelo type coincides with the regular type in the non integrable  case.
\begin{prop}\label{proptyp}
We have 
$$\Delta^1_{{\rm reg}}\left(\partial D,0\right) =\Delta^1\left(\partial D,0\right).$$ 
\end{prop}
\begin{proof}[Proof.]
We suppose that the origin is a point of finite  D'Angelo type. According to (\ref{eqcomp}) 
%$ \displaystyle \Delta^1_{{\rm reg}}\left(\partial D,0\right) \leq \Delta^1\left(\partial D,0\right)$ 
we may write: 
$$\Delta^1_{{\rm reg}}\left(\partial D,0\right)=2m<+\infty.$$ 
So we may assume that $u\left(\zeta\right)=\left(\zeta,0\right)$ is a regular 
$J$-holomorphic disc of maximal contact order $2m$, and that the structure $J$ satisfies (\ref{eqstr}) and (\ref{eqstrO}). 
Moreover the defining function $\rho$ has the following local expression: 
$$\rho=\Re ez_2+H_{2m}\left(z_1,\overline{z_1}\right)+O\left(|z_1|^{2m+1}+|z_2|\|z\|\right).$$

Now consider a $J$-holomorphic disc $v=\left(f_1,g_1,f_2,g_2\right): \left(\Delta,0\right) \rightarrow \left(\R^{4},0,J\right)$ of finite contact order 
satisfying $v\left(0\right)=0$ and such that $\delta\left(v\right)\geq 2$ (see definition \ref{defco} for notations). 

We set $v_1:=f_1+ig_1$ and $v_2:=f_2+ig_2$. The $J$-holomorphy equation for the disc $v$ is given by:
$$  
\left\{
\begin{array}{lll}  
\displaystyle a_k\left(v\right)\frac{\partial f_k}{\partial x}+b_k\left(v\right)\frac{\partial g_k}{\partial x}&= & 
\displaystyle\frac{\partial f_k}{\partial y},\\
& & \\
\displaystyle c_k\left(v\right)\frac{\partial f_k}{\partial x}-a_k\left(v\right)\frac{\partial g_k}{\partial x}&=& 
\displaystyle\frac{\partial g_k}{\partial y},\\
\end{array}
\right.
$$
for $k=1,2$. Since $J\left(v\right)=J_{st}+O\left(|v_2|\right)$ and $\delta\left(v\right)\geq 2$, it follows that: 
\begin{equation}\label{eqord}
\left\{
\begin{array}{lllll}  
\delta\left(v_1\right)&=&\delta\left(f_1\right)&=&\delta\left(g_1\right),\\
&&\\
 \delta\left(v_2\right)&=&\delta\left(f_2\right)&=&\delta\left(g_2\right).
\end{array}
\right.
\end{equation}
Then consider 
\begin{equation}\label{eqcon}
\rho\circ v\left(\zeta\right)=f_2\left(\zeta\right) +H_{2m}\left(v_1\left(\zeta\right),\overline{v_1\left(\zeta\right)}\right)+
O\left(|v_1\left(\zeta\right)|^{2m+1}+|v_2(\zeta)|\|v\left(\zeta\right)\|\right).
\end{equation}
Equation (\ref{eqord}) implies that the  term $O\left(|v_2|\|v\|\right)$ in (\ref{eqcon}) vanishes to 
 order larger than $f_2$. 

\vspace{0.7cm}

{\bf Case 1:} $\delta(f_2) > \delta\left(H_{2m}\left(v_1,\overline{v_1}\right)\right)$.  
In that case  
$$\delta_0\left(\partial D,u\right)=\delta\left(H_{2m}\left(v_1,\overline{v_1}\right)\right)=2m\delta\left(v_1\right).$$
Thus we get:
$$\frac{\delta_0\left(\partial D,v\right)}{\delta\left(v\right)}=\frac{2m\delta\left(v_1\right)}{\delta\left(v_1\right)}=2m.$$

\vspace{0.7cm}

{\bf Case 2:} $\delta(f_2) \leq \delta\left(H_{2m}\left(v_1,\overline{v_1}\right)\right)$. We have two subcases.  

\vspace{0.2cm}

{\bf Subcase 2.1:} $f_2+H_{2m}\left(v_1,\overline{v_1}\right)\not\equiv 0$.  Thus  
$$\delta_0\left(\partial D,u\right)=\delta\left(\Re e v_2\right)=\delta\left(v_2\right),$$
and so
$$\frac{\delta_0\left(\partial D,v\right)}{\delta\left(v\right)}=\frac{\delta\left(v_2\right)}{\delta\left(v\right)}\leq \frac{\delta\left(H_{2m}\left(v_1,\overline{v_1}\right)\right)}{\delta\left(v\right)}=\frac{2m\delta\left(v_1\right)}{\delta\left(v\right)}.$$
This means that: 
\begin{eqnarray*} 
 \frac{\delta_0\left(\partial D,v\right)}{\delta\left(v\right)}&=& 1 \mbox{ } \mbox{ if } \mbox{ } 
\displaystyle \delta\left(v\right)=\delta\left(v_2\right)
\end{eqnarray*}
or 
\begin{eqnarray*}
\displaystyle \frac{\delta_0\left(\partial D,v\right)}{\delta\left(v\right)}&\leq & 2m \mbox{ } \mbox{ if } \mbox{ }   
\displaystyle \delta\left(v\right)=\delta\left(v_1\right). 
\end{eqnarray*}

\vspace{0.2cm}

{\bf Subcase 2.2:} $f_2+H_{2m}\left(v_1,\overline{v_1}\right)\equiv 0$.
Let  $w:\Delta \rightarrow \left(\R^4,J_{st}\right)$ be a standard holomorphic disc satisfying 
$w\left(0\right)=0$ and:  
$$\frac{\partial^{k} w}{\partial x^{k}}\left(0\right) =  
\frac{\partial^{k} v}{\partial x^{k}}\left(0\right),$$
for $k=1,\cdots, 2m\delta\left(v\right)$.
Since $\delta\left(v_2\right)=2m\delta\left(v_1\right)=2m\delta\left(v\right)< +\infty$ and  
since $J\left(v\right)=J_{st}+O(|v_2|)$, any differentiation of $J\left(v\right)$, of order smaller than 
$2m\delta\left(v\right)$,  is equal to zero. 
Combining this with the $J$-holomorphy equation (\ref{eqholo}) of $v$
%$$J\left(v\right)\frac{\partial v}{\partial x}=\frac{\partial v}{\partial y},$$ 
we obtain: 
$$\frac{\partial^{k+l} w}{\partial x^{k}\partial y^l}\left(0\right) =  
\frac{\partial^{k+l} v}{\partial x^{k}\partial y^l}\left(0\right),$$
for $k+l=1,\cdots, 2m\delta\left(v\right)$.
Since $\rho \circ v$ vanishes to an order greater than $2m\delta\left(v\right)$ at 0 and since 
it involves only the $2m\delta\left(v\right)$-jet of $v$, it follows that 
$\rho \circ w$ vanishes to an  order greater than $2m\delta\left(v\right)$ at 0. 
Finally we have constructed a standard holomorphic disc $w$ such that 
\begin{equation*}
\left\{
\begin{array}{lll}  
\delta\left(w\right)&=&\delta\left(v\right),\\
&&\\
\delta_0\left(\partial D,w\right)&>&2m\delta\left(w\right),
\end{array}
\right.
\end{equation*}
 which is not possible since, according Proposition \ref{prop1},  the type 
for the standard structure of $\partial D$ at 
the origin is equal to $2m$.    
\end{proof}

\subsection{Construction of a local peak plurisubharmonic function}

We first give the definition  of a local peak $J$-plurisubharmonic function for a domain $D$. 
\begin{defi}
Let $D$ be a domain in an almost complex manifold $\left(M,J\right)$. A function $\varphi$ 
is called a local peak $J$-plurisubharmonic function at 
a boundary point $p \in \partial D$ if there exists a neighborhood $U$ of $p$ such that $\varphi$ is 
continuous up to $\overline{D}\cap U$ and satisfies: 
\begin{enumerate}
\item $\varphi$ is $J$-plurisubharmonic on  $D\cap U$,
\item $\varphi\left(p\right)=0$,
\item $\varphi<0$ on  $\overline{D}\cap U\backslash\{p\}.$
\end{enumerate} 
\end{defi}

%%\bigskip
%%\begin{center}
%%\input{peakf2.pstex_t}
%%\end{center}
%%\bigskip

The existence of  local peak $J_{st}$-plurisubharmonic functions was first proved by \\
E.Fornaess and N.Sibony in \cite{fo-si}.
For almost complex manifolds the existence was proved by S.Ivashkovich and J.-P.Rosay in \cite{iv-ro} 
whenever the domain is strictly $J$-pseudoconvex.
In the next Proposition we state the existence for $J$-pseudoconvex regions of finite D'Angelo type.  
As mentionned earlier our the considerations are purely local. In particular, the assumptions of $J$-plurisubharmonicity and of finite D'Angelo type may be 
restricted to a neighborhood of a boundary point. For convenience of writing, we state them globally.
\begin{theo}\label{proppeak}
Let $D=\{\rho<0\}$ be a domain of finite D'Angelo type in a four dimensional almost complex manifold $\left(M,J\right)$. 
We suppose that $\rho$ is a $\mathcal{C}^{2}$ defining function of $D$, $J$-plurisubharmonic on a neighborhood of $\overline{D}$. 
Let $p \in \partial D$ be a boundary point. Then  there exists a local peak $J$-plurisubharmonic function at $p$. 
\end{theo}
\begin{proof}[Proof.]
Since the existence of a local peak function near a boundary point of type $2$ was proved in \cite{iv-ro}, we assume that 
$p$ is a boundary point of D'Angelo type $2m>2.$ The problem being purely local we assume that $D\subset \C^2$ and 
that $p=0$. According to Proposition \ref{prop1} the defining 
function $\rho$ has the following local expression on a neighborhood $U$ of the origin: 
$$\rho=\Re ez_2+H_{2m}\left(z_1,\overline{z_1}\right)+\widetilde{H}(z_1,z_2)+O\left(|z_1|^{2m+1}+|z_2||z_1|^m+|z_2|^2\right)$$ 
where $H_{2m}$ is a subharmonic polynomial containing a nonharmonic part, denoted by  $H^{*}_{2m}$, and 
$$\widetilde{H}(z_1,z_2)=\Re e \displaystyle \sum_{k=1}^{m-1} \rho_{k} z_1^k z_2.$$
 
According to \cite{fo-si} (see Lemma 2.4), the polynomial $H_{2m}$ satisfies the following Lemma:    
\begin{lem}\label{lemfosi}
There exist a positive $\delta>0$ and a smooth function $g~: \R\rightarrow \R$ with period $2\pi$ with the following properties: 
\begin{enumerate}
\item $-2<g\left(\theta\right)<-1$,
\item $\|g\|<1/\delta$,
\item ${\rm max} \left(\Delta H_{2m},\Delta\left(\|H_{2m}^*\|g\left(\theta\right)|z_1|^{2m}\right)\right)>\delta\|H_{2m}^*\||z_1|^{2\left(m-1\right)}$, for $z_1=|z_1|e^{i\theta}\neq 0$ and,
\item $\Delta\left(H_{2m}+\delta\|H_{2m}^*\|g\left(\theta\right)|z_1|^{2m}\right)>\delta^2\|H_{2m}^*\||z_1|^{2\left(m-1\right)}$.
\end{enumerate}
\end{lem}

%*************** SIGNIFICATION GEOMETRIQUE~: FIGURE   ******************

We denote by $P$ the function defined by 
$$P\left(z_1,\overline{z_1}\right):=H_{2m}\left(z_1,\overline{z_1}\right)+\delta\|H^*_{2m}\|g\left(\theta\right)|z_1|^{2m}.$$
Theorem \ref{proppeak} will be proved by establishing the following claim. 
\vspace{0.7cm}\\
{\bf Claim.}
There are  positive constants $L$ and $C$ such that the function 
$$\varphi:=\Re e z_2 +2L\left(\Re e z_2\right)^2-L\left(\Im m z_2\right)^2+P(z_1,\overline{z_1})+\widetilde{H}(z_1,z_2)+C|z_1|^2|z_2|^2$$ 
is a local peak $J$-plurisubharmonic function at the origin. 
\vspace{0.7cm} \\
{\it Proof of the claim.}
We first prove that  the function $\varphi$ is $J$-plurisubharmonic. We set:
$$dd^{c}_{J}\varphi=\alpha_1 dx_1\wedge dy_1+\alpha_2 dx_2\wedge dy_2+\alpha_3 dx_1\wedge dx_2+\alpha_4 dx_1\wedge dy_2+
\alpha_5 dy_1\wedge dx_2+\alpha_6 dy_1\wedge dy_2,$$
where $\alpha_k$, for $k=1,\cdots,6,$ are real valued function. 
According to the matricial representation of $J$ (see (\ref{eqstr})), the Levi form of $\varphi$ at a point $z \in D\cap U$ 
and at a vector  $v=\left(X_1,Y_1,X_2,Y_2\right) \in T_z \R^4$ can be written 
\begin{eqnarray*}
\mathcal{L}_J \varphi\left(z,v\right)& = & c_1\alpha_1X_1^2-2a_1\alpha_1X_1Y_1-b_1\alpha_1Y_1^2+ \beta_3X_1X_2+\beta_4X_1Y_2+\\
&&\\
& & \beta_5Y_1X_2+ \beta_6Y_1Y_2+ c_2\alpha_2X_2^2-2a_2\alpha_2X_2Y_2-b_2\alpha_2Y_2^2,\\
\end{eqnarray*}
with
$$\left\{
\begin{array}{lll}
\beta_3 &~:= & \alpha_3\left(a_2-a_1\right)+\alpha_4c_2-\alpha_5c_1\\
\\
\beta_4 &~:= & -\alpha_4\left(a_1+a_2\right)+\alpha_3b_2-\alpha_6c_1\\
\\
\beta_5 &~:= & \alpha_5\left(a_1+a_2\right)-\alpha_3b_1+\alpha_6c_2\\
\\
\beta_6 &~:= & \alpha_6\left(a_1-a_2\right)-\alpha_4b_1+\alpha_5b_2.\\
\\
\end{array}\right.
$$
Moreover due to (\ref{eqstrO}) we have for  $k=1,2$
\begin{equation*}
\left\{
\begin{array}{lll}  
a_k=O\left(|z_2|\right)\\
\\
b_k=-1+O\left(|z_2|\right)\\
\\
c_k=1+O\left(|z_2|\right).
\end{array}
\right.
\end{equation*}
%$$a_k=O\left(|z_2|\right)\mbox {, } b_k=-1+O\left(|z_2|\right) \mbox { and } \mbox { }c_k=1+O\left(|z_2|\right) \mbox {, } k=1,2,$$ 
This implies that for $k=1,2$:
$$ c_k\alpha_kX_k^2-2a_k\alpha_kX_kY_k-b_k\alpha_kY_k^2\geq \frac{\alpha_k}{2}\left(X_k^2+Y_k^2\right).$$ 
Thus we obtain 
\begin{eqnarray*}
\mathcal{L}_J \varphi\left(z,v\right)& \geq & 
\frac{\alpha_1}{4}X_1^2+\beta_3X_1X_2+ \frac{\alpha_2}{4}X_2^2+
\frac{\alpha_1}{4}Y_1^2+  \beta_5Y_1X_2+\frac{\alpha_2}{4}X_2^2+\\
&&\\
&&\frac{\alpha_1}{4}X_1^2+\beta_4X_1Y_2+\frac{\alpha_2}{4}Y_2^2+ 
\frac{\alpha_1}{4}Y_1^2+ \beta_6Y_1Y_2+\frac{\alpha_2}{4}Y_2^2.
\end{eqnarray*}
In order to prove that  $\varphi$ is $J$-plurisubharmonic, we need to see that:
\begin{enumerate}
\item $\alpha_k\geq 0$,  for $k=1,2$,\\
\item $4\beta_j^2\leq \alpha_1\alpha_2$, for $j=3,\cdots,6$. 
\end{enumerate}  

%We have  noticed in the proof of Proposition \ref{prop1}  that since the Levi form for the standard 
%structure of $\widetilde{H}(z_1,z_2)$ is identically equal to zero, 
%the  levi form for $J$ of $\widetilde{H}(z_1,z_2)$ is obtained by either 
%differentiating two times $\widetilde{H}(z_1,z_2)$ with respect 
%$J-J_{st}=O(|z_2|)$ or differentiating $\widetilde{H}(z_1,z_2)$ and the structure. 

 The coefficient $\alpha_2$ is obtained by the differentiation of $\Re e z_2$, $2L\left(\Re e z_2\right)^2-L\left(\Im m z_2\right)^2$, 
$\widetilde{H}(z_1,z_2)$ and $C|z_1|^2|z_2|^2$. Hence we have for $z$ sufficiently close to the origin
\begin{eqnarray*}
\alpha_2 \geq L>0.
\end{eqnarray*} 

The coefficient $\alpha_1$ is obtained by differentiating $P$, $\widetilde{H}(z_1,z_2)$  and $C|z_1|^2|z_2|^2$. 
This is equal to 
\begin{eqnarray*}
\alpha_1 &=& \Delta P +O(|z_1|^{2m-2}|z_2|)+O(|z_2|^2)+C|z_2|^2+O(|z_2|^3)\\
&&\\
& \geq & \frac{\delta^2\|H_{2m}^*\|}{2}|z_1|^{2m-2}+ \frac{C}{2}|z_2|^2,
\end{eqnarray*} 
for $z$ sufficiently small and $C>0$ large enough. Hence $\alpha_1$ is nonnegative. 

Finally it sufficient to prove that 
\begin{equation*}
4\beta_j^2\leq   L\left(\frac{\delta^2\|H_{2m}^*\|}{2}|z_1|^{2m-2}+ \frac{C}{2}|z_2|^2\right),
\end{equation*}
to insure the $J$-plurisubharmonicity of $\varphi$.
The coefficient $|\beta_j|$ is equal to 
\begin{eqnarray*}
|\beta_j| &=& O(|z_2|)+LO(|z_2|^2)+O(|z_1|^{2m-1})+CO(|z_1||z_2|)\\
&&\\
& \leq & C'(|z_2|+|z_1|^{2m-1}),
\end{eqnarray*} 
for a positive constant $C'$ (not depending on  $L$ and $C$). It follows that $ \varphi$ is $J$-plurisubhar\-monic 
on a neighborhood of the origin.

\vspace{0.7cm}

We prove now that $\varphi$ is local peak at the origin, that is there exists $r>0$ such that $\overline{D} \cap \{0<\|z\|\leq r\}\subset \{\varphi<0\}$. 
Assuming that   $z \in \{\rho=0\}\cap  \{0<\|z\|\leq r\}$ we have:             
\begin{eqnarray*}
\varphi\left(z\right)&=&\delta\|H^*_{2m}\|g(\theta)|z_1|^{2m}+2L\left(\Re e z_2\right)^2-
L\left(\Im m z_2\right)^2+C|z_1|^2|z_2|^2+\\
&&\\
&&O\left(|z_1|^{2m+1}\right)+O\left(|z_2||z_1|^m\right)+O\left(|z_2|^2\right).
\end{eqnarray*}
Since $g<-1$ and increasing $L$ if necessary we have 
$$O\left(|\Im m z_2||z_1|^{m}\right) \leq -\frac{1}{2}\delta\|H^*_{2m}\|g\left(\theta\right)|z_1|^{2m}+\frac{1}{2}
L\left(\Im m z_2\right)^2,$$
whenever $z$ is sufficiently close to the origin. Thus
\begin{eqnarray*}
\varphi\left(z\right) &\leq & -\frac{1}{2}\delta\|H^*_{2m}\||z_1|^{2m}+(2L+C|z_1|^2)\left(\Re ez_2\right)^2-\frac{1}{2}L\left(\Im mz_2\right)^2+
C|z_1|^2(\Im m z_2)^2+\\
&&\\
& & O\left(|z_1|^{2m+1}\right)+O\left(|\Re e  z_2||z_1|^m\right)+O\left(|z_2|^2\right)\\
&&\\
&\leq &-\frac{1}{4}\delta\|H^*_{2m}\||z_1|^{2m}+(2L+C|z_1|^2)\left(\Re ez_2\right)^2-\frac{1}{4}L\left(\Im mz_2\right)^2+
O\left(|\Re e  z_2||z_1|^m\right)+\\
&&\\
& & O\left(|z_2|^2\right).
\end{eqnarray*}
There is a positive constant $C''$ such that  
$$O\left(|z_2|^2\right)\leq  C''|\Re e z_2|^2+ C'' |\Im m z_2|^2.$$
Thus  increasing $L$ if necessary:
\begin{eqnarray*}
\varphi\left(z\right) &\leq & 
-\frac{1}{4}\delta\|H^*_{2m}\||z_1|^{2m}+(2L+C|z_1|^2)\left(\Re ez_2\right)^2+O (|\Re e z_2|^2)\\
&&\\
& & - \left(\frac{1}{4}L-C''\right)\left(\Im m z_2\right)^2+
O\left(|\Re e  z_2||z_1|^m\right)+O(|\Im  mz_2|^2\|z\|).\\
&&\\
&\leq & -\frac{1}{4}\delta\|H^*_{2m}\||z_1|^{2m}+(2L+C|z_1|^2)\left(\Re ez_2\right)^2+O (|\Re e z_2|^2) +
O\left(|\Re e  z_2||z_1|^m\right)\\
&&\\
& & -\frac{1}{2}\left(\frac{1}{4}L-C''\right)\left(\Im m z_2\right)^2.
\end{eqnarray*}
Since
$$-\Re ez_2(1+O(|z|))=H_{2m}\left(z_1,\overline{z_1}\right)+O\left(|z_1|^{2m+1}+|\Im mz_2||z_1|+|\Im m z_2|^2\right),$$ 
we have
$$(\Re ez_2)^2(1+O(|z|)) = O\left(|z_1|^{4m}+|\Im m z_2||z_1|^{2m+1}+|\Im m z_2|^2\|z\|\right).$$ 
We finally obtain for $z$ small enough
\begin{eqnarray*}
\varphi\left(z\right)\leq -\frac{1}{8}\delta\|H^*_{2m}\||z_1|^{2m}-\frac{1}{4}\left(\frac{1}{4}L-C''\right)\left(\Im mz_2\right)^2.
\end{eqnarray*}
Thus $\varphi$ is negative for $z \in \{\rho=0\}\cap  \{0<\|z\|\leq r\}$, with $r$ small enough. 
It follows that, reducing $r$ if necessary, 
$$\overline{D} \cap \{0<\|z\|\leq r\}\subset \{\varphi<0\},$$ 
which achieves the proof of the claim and of Theorem \ref{proppeak}.
\end{proof}

We notice that in case  $\mathcal{L}_J\Re e z_2\equiv 0$, we may give  a simpler expression for a 
local peak $J$-plurisubharmonic function. 
\begin{prop}
If  $\mathcal{L}_J\Re e z_2\equiv 0$, then there exists  a real positive number $L$ such that the function 
$$\varphi:=\Re e z_2+2L\left(\Re e z_2\right)^2-L\left(z_2\right)^2+P\left(z_1,\overline{z_1}\right)$$ 
is local peak $J$-plurisubharmonic at the origin.
\end{prop}

\section{Estimates of the Kobayashi pseudometric}
In this section we prove standard estimates of the Kobayashi pseudometric on 
$J$-pseudocon\-vex regions  of finite D'Angelo type in an almost complex manifold.
\subsection{The Kobayashi pseudometric}
The existence of local pseudoholomorphic discs proved in \cite{ni-wo} 
allows to define the {\it Kobayashi pseudometric} $K_{\left(M,J\right)}$ for $p\in M$ and $v \in T_pM$ :
$$K_{\left(M,J\right)}\left(p,v\right):=\inf 
\left\{\frac{1}{r}>0, u : \Delta \rightarrow \left(M,J\right) 
\mbox{  $J$-holomorphic }, u\left(0\right)=p, d_{0}u\left(\partx\right)=rv\right\}.$$

Since the composition of pseudoholomorphic maps is still pseudoholomorphic, the 
 Kobayashi (infinitesimal) pseudometric satisfies the following decreasing property : 
\begin{prop}\label{propdec}
Let $f : \left(M',J'\right)\rightarrow \left(M,J\right)$ be a $\left(J',J\right)$-holomorphic map. Then for any
 $p \in M'$ and $v \in T_{p}M'$ we have 
$$K_{\left(M',J'\right)}\left(p,v\right)\geq K_{\left(M,J\right)}\left(f\left(p\right),d_{p}f\left(v\right)\right).$$   
\end{prop}

Let $d_{\left(M,J\right)}$ be the  integrated pseudodistance of $K_{\left(M,J\right)}$ :
$$d_{\left(M,J\right)}\left(p,q\right): =\inf\left\{\int_0^1 
K_{\left(M,J\right)}\left(\gamma\left(t\right),\dot{\gamma}\left(t\right)\right)dt, \mbox{ }
\gamma : [0,1]\rightarrow M, \mbox{ }\gamma\left(0\right)=p, \gamma\left(1\right)=q\right\}.$$

Similarly to the standard  integrable case, B.Kruglikov (see \cite{kr2}) proved  that
the integrated pseudodistance of the Kobayashi pseudometric coincides with the  Kobayashi pseudodistance 
defined by chains of pseudholomorphic discs.  
Then we define :  

\begin{defi}\label{defin}
\begin{enumerate}
\item The manifold $\left(M,J\right)$ is Kobayashi hyperbolic if the integrated pseudodistance $d_{\left(M,J\right)}$ is a distance.
\item The manifold $\left(M,J\right)$ is local Kobayashi hyperbolic at $p \in M$ if there exist a neighborhood $U$ of $p$ and a positive constant $C$ such that 
$$K_{\left(M,J\right)}\left(q,v\right)\geq C\|v\|$$ for every $q \in U$ and every $v \in T_qM$.  
\item A Kobayashi hyperbolic manifold $\left(M,J\right)$ is complete hyperbolic if it is complete for the distance  $d_{\left(M,J\right)}$.
\end{enumerate}
\end{defi}
\subsection{Hyperbolicity of pseudoconvex regions of finite D'Angelo type}

In order to localize pseudoholomorphic discs, we need the following technical Lemma (see \cite{ga-su} for a proof).  
\begin{lem}\label{lemlem}
Let $0 < r < 1$ and let $\theta_r$ be a smooth nondecreasing function on
$\R^+$ such that $\theta_r\left(s\right)= s$ for $s \leq r/3$ and $\theta_r\left(s\right) = 1$ 
for $s \geq 2r/3$. Let $\left(M,J\right)$ be an almost complex manifold, and let $p$ be a point of $M$.
Then there exist a neighborhood $U$ of $p$, 
positive constants $A = A\left(r\right)\geq 1$, $B=B\left(r\right)$, and a diffeomorphism $z~: U \rightarrow \B$ 
such that $z\left(p\right)=0$, $z_{*}J\left(p\right)=J_{st}$ and the function ${\rm log}\left(\theta_r\left(\vert z
\vert^2\right)\right) + \theta_r\left(A\vert z \vert\right) + B\vert z \vert^2$ is
$J$-plurisubharmonic on $U$.
\end{lem}

In the next Proposition we give a priori estimates and a localization principle  of the Kobayashi pseudometric. This
proves the local Kobayashi hyperbolicity of $J$-pseudoconvex $\mathcal{C}^2$  regions of finite D'Angelo type. 
If $(M,J)$ admits a global $J$-plurisubharmonic function, then K.Diederich and  A.Sukhov proved in \cite{di-su}
the (global) Kobayashi hyperbolicity of a relatively compact $J$-pseudoconvex domains (with $\mathcal{C}^3$ boundary)
 by constructing a bounded strictly $J$-plurisubharmonic exhaustion function. We notice that, in our case, if the manifold 
$(M,J)$ admits a global $J$-plurisubharmonic function then $J$-pseudoconvex $\mathcal{C}^2$ relatively compact regions 
of finite D'Angelo type are also (globally)  Kobayashi hyperbolic. 

\begin{prop}\label{thm}
Let $D=\{\rho<0\}$ be a  domain of finite D'Angelo type in an almost complex manifold $\left(M,J\right)$, 
where $\rho$ is a $\mathcal{C}^{2}$ defining function of $D$, $J$-plurisubharmonic 
in a neighborhood of $\overline{D}$. Let $p \in \bar{D}$ and  let $U$ be a neighborhood of $p$ in $M$. 
Then there exist   positive constants $C$ and $s$, and a neighborhood $V \subset U$ of $p$ 
in $M$, such that for each $q \in D \cap V$ and each $v\in T_q M$: 
\begin{equation}\label{ee1}
 K_{\left(D,J\right)}\left(q,v\right)\geq C \|v\|,
\end{equation}

\begin{equation}\label{ee2}
 K_{\left(D,J\right)}\left(q,v\right)\geq s K_{\left(D \cap U,J\right)}\left(q,v\right).
\end{equation}
\end{prop}
This Proposition is a classical application of Lemma \ref{lemlem}. This is due to N.Sibony \cite{nsi} 
(see also  \cite{be1}  and \cite{ga-su} for a proof). For convenience
we give the proof. 
\begin{proof}[Proof.]
According to Theorem \ref{proppeak}, there exists a local peak $J$-plurisubharmonic 
function $\varphi$  at $p$ for $D$. 
We can choose constants $0<\alpha<\alpha'<\beta'<\beta$ and $N>0$ such that $\varphi\geq -\beta^2/N$ on 
$\{ \|z\|<\alpha\}$ and 
$\varphi \leq -2\beta^2/N$ on $\overline{D}\cap \{\alpha' \leq \|z\|\leq \beta'\}$.

We define $\tilde{\varphi}$ by: 
$$
 \tilde{\varphi}:=
\left\{
\begin{array}{ll}  
{\rm max}\left(N\varphi+\|z\|^{2}-\beta^2,-2\beta^2\right) & \mbox{ if } z \in D\cap \{\|z\|\leq \beta'\},\\
&  \\
 -2\beta^2 & \mbox{ on }  D \backslash \{\|z\|\leq \beta'\}.\\
\end{array}\right.
$$
The function $\|z\|^{2}$ is $J$-plurisubharmonic on $\{q \in U~: |z\left(q\right)|<1\}$ if $\|z_{*}J-J_{st}\|_{\mathcal{C}^{2}\left(\B\right)}$ is 
sufficiently small. Then it follows that $\tilde{\varphi}$ is $J$-plurisubharmonic on $D.$ 
We may also suppose that $\tilde{\varphi}$ is negative on $D.$ 
Moreover the function $\tilde{\varphi}-\|z\|^{2}$ is $J$-plurisubharmonic on $D\cap \{q \in U~: |z\left(q\right)|\leq \alpha\}.$

Let $\theta_{\alpha^2}$ be a smooth non decreasing function on $\R^+$ such that $\theta_{\alpha^2}\left(s\right)=s$ for 
$s\leq \alpha^2/3$ and $\theta_{\alpha^2}\left(s\right)=1$ for $s\geq 2\alpha^2/3$. 
Set $V=\{q \in U~: |z\left(q\right)|\leq \alpha^2\}$. According to Lemma \ref{lemlem}, 
there are uniform positive constants $A\geq 1$ and $B$ such that the function 
$$
{\rm log}\left(\theta_{\alpha^2}\left(|z-z\left(q\right)|^2\right)\right)+ \theta_{\alpha^2}\left(A|z-z\left(q\right)|\right)+ B\|z\|^2
$$
is $J$-plurisubharmonic on $U$ for every $q \in D\cap V$. 

We define for each 
$q \in D\cap V$ the function:   

$$
\Psi_{q}:=
\left\{
\begin{array}{ll}  
 \theta_{\alpha^2}\left(|z-z\left(q\right)|^2\right){\rm exp}\left(\theta_{\alpha^2}\left(A|z-z\left(q\right)|\right)\right) 
{\rm exp}\left(B \tilde{\varphi}\left(z\right)\right) & \mbox{ on } D\cap \{\|z\|<\alpha\},\\
 & \\
\exp \left(1+B\tilde{\varphi}\right) & \mbox{ on } D\setminus \{\|z\|<\alpha\}.\\ 
\end{array}
\right.
$$
The function ${\rm log}\Psi_{q}$ is $J$-plurisubharmonic on $D\cap \{\|z\|<\alpha\}$ and, on   
 $D\setminus \{\|z\|<\alpha\}$, it coincides with $1+B\tilde{\varphi}$ which is 
 $J$-plurisubharmonic. Finally ${\rm log} \Psi_{q}$ is $J$-plurisubharmonic on the whole domain $D$.

Let $q \in V$ and let $v \in T_q M$ and consider a $J$-holomorphic disc $u~: \Delta \rightarrow D$  
such that $u\left(0\right)=q$ and
$d_0u \left(\partial /\partial x\right)= rv$ where $r>0$. For $\zeta$ sufficiently close to 0 we have
$$
u\left(\zeta\right) = q + d_0u\left(\zeta\right) +
\mathcal O\left(|\zeta|^2\right).
$$
We define the following function
$$
\phi\left(\zeta\right):= \frac{\Psi_q\left(u\left(\zeta\right)\right)}{|\zeta|^2}
$$
which is subharmonic on
$\Delta \backslash \{0\}$ since ${\rm log} \phi$ is subharmonic. 
If $\zeta$ close to $0$, then   
\begin{equation}\label{eq0}
\phi\left(\zeta\right) = \frac{|u\left(\zeta\right)-q|^2}{|\zeta|^2} \exp\left(A|u\left(\zeta\right)-q|\right) 
\exp\left(B \tilde{\varphi}\left(u\left(\zeta\right)\right)\right).
\end{equation}
Setting $\zeta= \zeta_1+i\zeta_2$ and using
the $J$-holomorphy condition $d_0u\circ J_{st} = J \circ
d_0u$, we may write~:
$$
d_0u\left(\zeta\right) = \zeta_1 d_0u\left(\partial /
\partial x\right) + \zeta_2 J\left(d_0u\left(\partial / \partial x\right)\right).
$$
\begin{equation}\label{eqdis}
|d_0u\left(\zeta\right)| \leq |\zeta| \left(\|I+J\|\,\|d_0u\left(\partial
/\partial x\right)\|\right)
\end{equation}
According to (\ref{eq0}) and to (\ref{eqdis}), we obtain that $\limsup_{\zeta \rightarrow 0}\phi\left(\zeta\right)
$ is finite. Moreover setting $\zeta_2=0$ we have 
$$
\limsup_{\zeta \rightarrow 0}\phi\left(\zeta\right) \geq \|d_0u\left(\partial
/\partial x\right)\|^2\exp\left(B\tilde{\varphi}\left(q\right)\right).
$$
Applying the maximum principle to a subharmonic extension of $\phi$
on $\Delta$ we obtain the inequality 
$$
\|d_0u\left(\partial / \partial x\right)\|^2 \leq \exp\left(1-B\tilde{\varphi}\left(q\right)\right).
$$

Hence, by definition of the Kobayashi  pseudometric, 
we obtain for every $q \in D \cap V$ and every $v \in T_q M$:
\begin{equation*}
K_{\left(D,J\right)}\left(q,v\right) \geq  \left(\exp\left(-1+B\tilde{\varphi}\left(q\right)\right)\right)^{\frac{1}{2}}\|v\|.
\end{equation*}
This gives estimate (\ref{ee1}).

\vspace{0.5cm}

Now in order to obtain estimate (\ref{ee2}), we prove that there is a neighborhood $V\subset U$ and a 
positive constant $s$ such that for any $J$-holomorphic disc $u~: \Delta \rightarrow D$ with $u\left(0\right) \in V$ then 
$u\left(\Delta_s\right)\subset D\cap U.$ Suppose this is not the case. We obtain a sequence $\zeta_{\nu}$ of $\Delta$  and a 
sequence of $J$-holomorphic discs $u_{\nu}$ such that 
 $\zeta_{\nu}$ converges to 0, $u_{\nu}\left(0\right)$ converges to $p$ 
and $\|u_{\nu}\left(\zeta_{\nu}\right)\| \notin D\cap U$ for every $\nu.$ 
According to the estimate (\ref{ee1}), we obtain for a positive constant $c>0$:  
$$c \leq  d_{\left(D,J\right)}\left(u_{\nu}\left(0\right),u_{\nu}\left(\zeta_{\nu}\right)\right) 
\leq d_{\Delta}\left(\zeta_{\nu},0\right).$$
This contradicts the fact that  $\zeta_{\nu}$ converges to 0.
\end{proof}
 
The (global)  Kobayahsi hyperbolicity is provided if we suppose  that there is a global strictly $J$-plurisubharmonic function on 
$\left(M,J\right)$. 
\begin{cor}\label{2corhyphyp}
Let $D=\{\rho<0\}$ be a relatively compact domain of finite D'Angelo type in an almost 
complex manifold $\left(M,J\right)$ of dimension four, $\rho$ being a defining function of $D$, $J$-plurisubharmonic 
in a neighborhood of $\overline{D}$. Assume that $\left(M,J\right)$ admits a global strictly $J$-plurisubharmonic function. 
Then $\left(D,J\right)$ is Kobayahsi hyperbolic.
\end{cor}
As an application of the a priori estimate (\ref {ee1}) of Proposition \ref{thm}, we prove the tautness of $D$.
\begin{cor}
Let $D=\{\rho<0\}$ be a relatively compact  domain of finite D'Angelo type in an almost 
complex manifold $(M,J)$ of dimension two. Assume that $\rho$ is 
$J$-plurisubhar\-monic in a neighborhood of $\overline{D}$. Moreover suppose that 
$\left(M,J\right)$ admits a global strictly $J$-plurisubharmonic function. 
Then $D$ is taut.
\end{cor}
\begin{proof}[Proof.]
Let $(u_\nu)_\nu$ be a sequence of $J$-holomorphic discs in $D$. According to Corollary \ref{2corhyphyp} the domain 
$D$ is hyperbolic. Thus  the sequence $(u_\nu)_\nu$ is equiconituous, and then by Ascoli Theorem, 
we can extract from this sequence a subsequence still denoted  $(u_\nu)_\nu$ which converges to a map 
$u~: \Delta\rightarrow \overline{D}$. Passing to the limit the equation of 
$J$-holomorphicity of each $u_\nu$, it follows that $u$ is a $J$-holomorphic disc. Since  $\rho$ is 
$J$-plurisubharmonic defining function for $D$, we have, by applying the maximun principle to $\rho \circ u$, 
the alternative: either $u(\Delta)\subset D$ or  $u(\Delta)\subset \partial D$. 
\end{proof}

We point out that the tautness of the domain $D$ was proved, using a diferent method, by 
K.Diederich-A.Sukhov in \cite{di-su}.

\subsection{Uniform estimates of the Kobayashi pseudometric}

In order to obtain more precise estimates, we need to uniform estimates (\ref{ee1}) of the 
Kobaya\-shi pseudometric for a sequence of domains. 

\begin{prop}\label{propstab}
Assume that $D=\{\Re e z_2+P\left(z_1,\overline{z_1}\right)<0\}$ is a $J_{st}$-pseudoconvex region of $\R^4$, 
where $P$ is a homogeneous polynomial of degree $2k\leq 2m$ admitting a nonharmonic part.  
Let $D_{\nu}$ be a  sequence of $J_{\nu}$-pseudoconvex region of $\R^4$ such that $0 \in \partial D_{\nu}$ 
is a boundary point of finite D'Angelo type $2l_{\nu}\leq 2m$.
Suppose that $D_{\nu}$ converges in the sense of local Hausdorff set convergence to 
$D$ when $\nu$ tends to $+\infty$ and that $J_{\nu}$ converges to $J_{st}$ in the 
$\mathcal{C}^2$ topology when $\nu$ tends to $+\infty$.   
Then there exist a positive constant $C$ and a neighborhood $V \subset U$ of the origin 
 in $\R^4$, such that for large $\nu$ and for  every $q \in D_{\nu} \cap V$ and every $v\in T_q \R^4$  
\begin{equation*}
 K_{\left(D_{\nu},J\right)}\left(q,v\right)\geq C \|v\|.
\end{equation*}
\end{prop}

\begin{proof}[Proof.]
Under the conditions of Proposition \ref{propstab} we have the following Lemma: 
\begin{lem}
For every large $\nu$, there exists a diffeomorphism $\Phi_{\nu}~: \R^4 \rightarrow \R^4$ with the following property: 
\begin{enumerate} 
\item The map $\zeta\mapsto \left(\zeta,0\right)$ is a  $\left(\Phi_{\nu}\right)_{*}J_{\nu}$-holomorphic disc of 
maximal contact order $2l_{\nu}$. 
\item The almost complex structure $\left(\Phi_{\nu}\right)_{*}J_{\nu}$ satisfies conditions (\ref{eqstr}) and (\ref{eqstrO}).
\item $\Phi_{\nu}\left(D_{\nu}\right)=\{\rho_{\nu}<0\}$ with 
$$\rho_{\nu}=\Re e z_2+ \sum_{j=2l_{\nu}}^{2m}P_{j,\nu}\left(z_1,\overline{z_1}\right) 
+O\left(|z_1|^{2m+1}+|z_2|\|z\|\right)<0,$$ where  $P_{j,\nu}$ are homogeneous polynomials of degree $j$ and  $P_{2l_\nu,\nu}$ 
contains  a nonharmonic part denoted by   $P_{2l_\nu,\nu}^*\neq 0$.
\item we have $inf_{\nu}\{\|P_{2l_\nu,\nu}\|\}>0.$
\end{enumerate}
Moreover the sequence of diffeomorphisms $\Phi_{\nu}$ converges to the identity 
on any compact subsets of $\R^4$ in the $\mathcal{C}^2$ topology.
\end{lem}
The crucial fact used to prove Proposition \ref{propstab} is the point $\left(4\right)$, which is a direct consequence of the 
convergence of $\Phi_{\nu}\left(D_{\nu}\right)$ to $D$. Hence the proof of Proposition \ref{propstab} is similar to 
Theorem \ref{proppeak} and Theorem \ref{thm}, where all the constants are uniform.  

%This provides constants $L>0$, $\delta>0$ and functions $g_{\nu}$ satisfying Lemma \ref{lemfosi} 
%applied to $P_{2l_{\nu},\nu}^*$  such that the function defined by 
%$$\Psi_{\nu}\left(z\right)=\Re e z_2+2L\left(\Re e z_2\right)^2-L\left(z_2\right)^2+P_{2l_\nu,\nu}\left(z_1,\overline{z_1}\right)
%+\delta\|P^*_{2l_{\nu},\nu}\|g_{\nu}\left(\theta\right)|z_1|^{2l_{\nu}}+|z_2|^2|z_1|^2$$
%is a peak  $\left(\Phi_{\nu}\right)_{*}J_{\nu}$-plurisubharmonic function at the origin for $\Phi_{\nu}\left(D_{\nu}\right)$. Thus 
%$\Psi_{\nu}\circ \Phi_{\nu}$ is a peak  $J_{\nu}$-plurisubharmonic function at the origin for $D_{\nu}$.
%There exists constants $0<\alpha<\alpha'<\beta'<\beta$ and $N>0$ 
%such that $\Psi_{\nu}\geq -\beta^2/N$ on $\|z\|<\alpha$ and 
%$\Psi_{\nu} \leq -2\beta^2/N$ on $\overline{D}\cap \{\alpha' \leq \|z\|\leq \beta'\}$. 
%the rest of the proof is similar to  the Proof of Theorem \ref{thm} and gives the desired estimate.
\end{proof}

\subsection{H\"older extension of diffeomorphisms}

This subsection is devoted to the boundary continuity of diffeomorphisms. This is stated as follows:
\begin{prop}\label{corext}
Let $D=\{\rho<0\}$ and $D'=\{\rho'<0\}$ be two relatively compact domains of finite D'Angelo type $2m$ in four 
dimensional almost complex manifolds $\left(M,J\right)$ and $\left(M',J'\right)$. 
We suppose that $\rho$ (resp.  $\rho'$) is a $J$(resp $J'$)-plurisubharmonic defining function on a neighborhood 
of $D$  (resp. $\overline{D'}$). Let $f: D \rightarrow D'$ be a $\left(J,J'\right)$-biholomorphism. Then $f$ extends as a 
H\"older homeomorphism with exponent $1/2m$ between $\overline{D}$ and $\overline{D'}$.  
\end{prop}

Estimates of the Kobayashi pseudometric obtained by H.Gaussier and A.Sukhov 
in \cite{ga-su} provide the H\"older extension  with 
exponent $1/2$ up to the boundary of a biholomorphism between two strictly pseudoconvex 
domains (see Proposition 3.3 of \cite{co-ga-su}). Similarly, in order to obtain Proposition \ref{corext}, 
we begin by establishing a more precise estimate than (\ref{ee1}) of Proposition \ref{thm}.  
\begin{prop}\label{thmext}
Let $D=\{\rho<0\}$ be a  domain of finite D'Angelo type in a four dimensional almost complex manifold $\left(M,J\right)$, 
where $\rho$ is a $\mathcal{C}^{2}$ defining function of $D$, $J$-plurisubharmonic 
in a neighborhood of $\overline{D}$. Let $p \in \partial D$ and  let $U$ be a neighborhood of $p$ in $M$. 
Then there are   positive constant $C$ and a neighborhood $V \subset U$ of $p$ 
in $M$, such that for every $q \in D \cap V$ and every $v\in T_q M$: 
\begin{equation}\label{ee02}
 K_{\left(D,J\right)}\left(q,v\right)\geq C \frac{\|v\|}{{\rm dist}\left(q,\partial D\right)^{1/2m}}.
\end{equation}
\end{prop}

\begin{proof}[Proof of Proposition \ref{thmext}.]
  Let $p \in \partial D$. We may suppose that $D\subset \R^4$,  $p=0$ and that $J$ satisfies (\ref{eqstr}) and (\ref{eqstrO}). 
Let $q'$ be a boundary point in a neighborhood of 
the origin and let $\varphi_{q'}$ be the local peak $J$-plurisubharmonic function at $q'$ given by Theorem \ref{proppeak}.
There are positive constants $C_1$ and $C_2$ such that  
\begin{equation}\label{eqine}
-C_1\|z-q'\|\leq \varphi_{q'}\left(z\right) \leq -C_2\Psi_{q'}\left(z\right),
\end{equation}
where 
$$\Psi_{q'}\left(z\right):=|z_1-q_1'|^{2m}+|z_2-q_2'|^2+|z_1-q_1'|^2|z_2-q_2'|^2$$ 
is a $J$-plurisubharmonic function on a neighborhood $U$ of the origin.

Now consider a $J$-holomorphic disc $u~: \Delta \rightarrow D$,  such that $u\left(0\right)$ is sufficiently 
close to the origin and then, according to Proposition \ref{thm}, we have $\displaystyle u\left(\Delta_s\right)\subset D\cap U,$ 
for some $0<s<1$ depending only on $u\left(0\right)$. 
We assume that $q'$ is such that ${\rm dist}\left(u\left(0\right),\partial D\right)=\|u\left(0\right)-q'\|.$ 
According to the   $J$-plurisubharmonicity of $\Psi_{q'}$,  we have  for $|\zeta|\leq s$:
$$\Psi_{q'}\left(u\left(\zeta\right)\right)\leq \frac{C_3}{2\pi}\displaystyle 
\int_0^{2\pi}\Psi_{q'}\left(u\left(re^{i\theta}\right)\right)d\theta,$$
for some positive constant $C_3$.
Hence using (\ref{eqine}) and the $J$-plurisubharmonicity of $\varphi_{q'}$ we obtain:
$$\Psi_{q'}\left(u\left(\zeta\right)\right)\leq -\frac{C_3}{2\pi C_2}\displaystyle \int_0^{2\pi} \varphi_{q'}\left(u\left(re^{i\theta}\right)\right)d\theta
\leq -\frac{C_3}{C_2}\varphi_{q'}\left(u\left(0\right)\right).$$
Since there is a positive constant $C_4$ such that 
$$\|u\left(\zeta\right)-q'\|^{2m}\leq C_4\Psi_{q'}\left(u\left(\zeta\right)\right)$$ 
and using (\ref{eqine}), we finally obtain: 
$$\|u\left(\zeta\right)-q'\|^{2m}\leq \frac{C_1C_3C_4}{C_2}{\rm dist}\left(u\left(0\right),\partial D\right).$$

Hence there exists a positive constant $C_5$  such that: 
$${\rm dist}\left(u\left(\zeta\right),\partial D\right)\leq C_5 {\rm dist}\left(u\left(0\right),\partial D\right)^{1/2m},$$
whenever $\zeta\leq s$.

According to Lemma $1.5$ of \cite{iv-ro} there is a positive constant  $C_6$ such that:
$$\|\nabla u\left(0\right)\|\leq C_6 \sup_{|\zeta|<s} \|u\left(\zeta\right)-u\left(0\right)\|\leq  
C_5C_6 {\rm dist}\left(u\left(0\right),\partial D\right)^{1/2m},$$ 
 wich provides the desired estimate.
\end{proof}

We also need the two next lemmas provided by \cite{co-ga-su}:
\begin{lem}
\label{lemlowest1}
Let $D$ be a domain in  an almost complex manifold $\left(M,J\right)$. 
Then there is a positive constant $C$ such that for any $p \in D$ and any $v \in T_pM$:
\begin{equation}\label{eqlemlow}
K_{\left(D,J\right)}\left(p,v\right)\leq C \frac{\|v\|}{{\rm dist}\left(p,\partial D\right)}.
\end{equation}
\end{lem}

\vspace{0,3cm}

\begin{lem}\label{lemhopf}
(Hopf lemma) Let $D$ be a relatively compact domain with a $\mathcal C^2$
boundary on an almost complex manifold $(M,J)$. Then for any negative
$J$-plurisubharmonic function $\rho$ on $D$ there exists a constant $C > 0$ such that for any $p \in D$:
$$|\rho(p)| \geq C {\rm dist}(p,\partial D).$$ 
\end{lem}

\vspace{0,3cm}

Now we can go on the proof of Proposition \ref{corext}.
\begin{proof}[Proof of Proposition \ref{corext}]
Let  $f: D \rightarrow D'$ be a $\left(J,J'\right)$-biholomorphism. According to Proposition \ref{thmext} 
and to the decreasing property of the Kobayashi pseudometric
there is a positive constant $C$ such that 
for every $p \in D$ sufficiently close to the boundary  and every $v \in T_pM$
$$C\frac{\|d_pf\left(v\right)\|}{{\rm dist}\left(f\left(p\right),\partial D'\right)^{\frac{1}{2m}}}
\leq K_{\left(D',J'\right)}\left(f\left(p\right),d_pf\left(v\right)\right) 
= K_{\left(D,J\right)}\left(p,v\right).$$
Due to Lemma \ref{lemlowest1} 
there exists a positive constant $C_1$ such that: 
$$K_{\left(D,J\right)}\left(p,v\right)\leq C_1 \frac{\|v\|}{{\rm dist}\left(p,\partial D\right)}.$$
This leads to: 
$$\|d_pf\left(v\right)\| \leq \frac{C_1}{C}  \frac{{\rm dist}\left(f\left(p\right),
\partial D'\right)^{\frac{1}{2m}}}{{\rm dist}\left(p,\partial D\right)}\|v\|.$$
Moreover the Hopf lemma \ref{lemhopf} for almost complex manifolds applied to  $\rho' \circ f$ and $\rho
\circ f^{-1}$  and the fact that $\rho$ and $\rho'$ are defining functions, 
provides the following boundary distance preserving property:
$$ \frac{1}{C_2}{\rm dist}\left(p,\partial D\right) \leq {\rm dist}\left(f\left(p\right),\partial D'\right)
 \leq C_2{\rm dist}\left(p,\partial D\right),$$
for some positive constat $C_2$. Finally this implies:
$$\|d_pf\left(v\right)\| \leq \frac{C_1C_2}{C} \frac{\|v\|}{{\rm dist}\left(p,\partial D\right)^{\frac{2m-1}{2m}}}.$$
This gives the desired statement.
\end{proof}

\section{Sharp estimates of the Kobayashi pseudometric}

In this section we give sharp lower estimates of the Kobayashi pseudometric in a pseudoconvex region near a boundary point of 
finite D'Angelo type less than or equal to four. This condition will appear necessary, in our proof, as explained in the appendix. 
Moreover  in order to give sharp estimates near a point of arbitrary finite  D'Angelo type,  
we are also interested in the nontangential behaviour  of  the Kobayashi pseudometric. 

\vspace{0.2cm}

The main result of this section is the following theorem (see also Theorem B):
\begin{theo}\label{theo}
Let $D=\{\rho<0\}$ be a relatively compact domain  of finite D'Angelo type less than or equal to four 
in an almost complex manifold $\left(M,J\right)$ of dimension four, where $\rho$ is a 
$\mathcal{C}^{2}$ defining function of $D$, $J$-plurisubharmonic on a neighborhood of $\overline{D}$. 
Then there exists a positive constant $C$ with the following property: 
for every $p \in D$ and every $v \in T_{p}M$ there is a diffeomophism, $\Phi_{p^*}$, in a neighborhood $U$ of $p$, 
such that: 
\begin{equation}\label{est1}
K_{\left(D,J\right)}\left(p,v\right)\geq C\left(\frac{|\left(d_p\Phi_{p^*}v\right)_1|}{\tau\left(p^*,|\rho\left(p\right)|\right)}+\frac{|\left(d_p\Phi_{p^*}v\right)_{2}|}{|\rho\left(p\right)|}\right),
\end{equation}  
where $\tau\left(p^*,|\rho\left(p\right)|\right)$ is defined  by (\ref{eqtau}). 

As a direct consequence we have:
\begin{equation}\label{est2}
K_{\left(D,J\right)}\left(p,v\right)\geq C'\left(\frac{|\left(d_p\Phi_{p^*}v\right)_1|}{|\rho\left(p\right)|^{\frac{1}{4}}}+\frac{|\left(d_p\Phi_{p^*}v\right)_{2}|}{|\rho\left(p\right)|}\right),
\end{equation}
for a positive constant $C'$.
\end{theo} 

In complex manifolds, D.Catlin \cite{ca} first obtained
such an  estimate, based on lower estimates of the Carath\'eodory pseudometric. F.Berteloot \cite{be2} gave a different proof based 
on a Bloch principle. Our proof wich is inspired by the proof of F.Berteloot 
is based on some scaling method. 

\subsection{The scaling method}
We consider here a pseudoconvex region $D=\{\rho<0\}$ of finite D'Angelo type $2m$ in $\R^4$, where $\rho$ has the following 
expression on a neighborhood $U$ of the origin:   
$$\rho\left(z_1,z_2\right)=\Re ez_2+H_{2m}\left(z_1,\overline{z_1}\right)+O\left(|z_1|^{2m+1}+|z_2|\|z\|\right).$$
where $H_{2m}$ is a homogeneous subharmonic polynomial of degree $2m$ admitting a nonharmonic part. 

Assume that $p_{\nu}$ is a sequence of points in $D\cap U$ converging to the origin. For each 
$p_{\nu}$ sufficiently close to $\partial D$, 
there exists a unique point $p_{\nu}^{*} \in \partial D \cap U$ such that 
$$p_{\nu}^*=p_{\nu}+\left(0,\delta_{\nu}\right),$$
 with $\delta_{\nu}>0$. Notice that for large $\nu$, the quantity  $\delta_{\nu}$ is 
equivalent to ${\rm dist}\left(p_{\nu},\partial D \cap U\right)$ and to $|\rho\left(p_{\nu}\right)|$.

\vspace{0.5cm}

We consider a  diffeomorphism $\Phi^{\nu}: \R^4\rightarrow \R^4$  satisfying:
\begin{enumerate}
\item $\Phi^{\nu}\left(p_{\nu}^{*}\right)=0$ and $\Phi^{\nu}\left(p_{\nu}\right)=\left(0,-\delta_{\nu}\right)$. 
\item $\Phi^{\nu}$ converges to $Id~:\R^{4} \rightarrow \R^{4}$ on any compact subset of $\R^{4}$ in the $\mathcal{C}^{2}$ 
sense.
\item  When we denote by $D^{\nu}:=\Phi^{\nu}\left(D\cap U\right)$ which admits the defining function is  
$\rho^{\nu}:=\rho\circ \left(\Phi^{\nu}\right)^{-1}$ and by $J^{\nu}:=\left(\Phi^{\nu}\right)_*J$, then 
 $\rho^{\nu}$ is given by: 
$$\rho^{\nu}\left(z_1,z_2\right)= \Re ez_2+\displaystyle \sum_{k=2l_{\nu}}^{2m}P_{k}\left(z_1,\overline{z_1},p_{\nu}^*\right)
+O\left(|z_1|^{2m+1}+|z_2|\|z\|\right),$$
where the polynomial $P_{2l_{\nu}}$  contains a nonharmonic part. Moreover $J^{\nu}$ satisfies (\ref{eqstr}) 
and (\ref{eqstrO}). 
\end{enumerate} 
This is done by considering first the translation $T^{\nu}$ of $\R^4$ given by $z\mapsto z-p_{\nu}^*$. 
According to J.-F.Barraud and E.Mazzilli \cite{ba-ma} that the D'Angelo type is an upper
semicontinuous function in a four dimensional almost complex manifold. 
Thus the D'Angelo type of points in a small enough neighborhood can only be smaller than at the point itself. 
Then we consider  a $\left(T^{\nu}\right)_*J$-holomorphic disc  $u$ of maximal contact order $2l_{\nu}$, 
where $2l_{\nu}\leq 2m$ is the D'Angelo type of $p_{\nu}^*$. We choose coordinates such that $u$ is 
given by $u\left(\zeta\right)=\left(\zeta,0\right)$, and such that $\left(T^{\nu}\right)_*J\left(z_1,0\right)=J_{st}$ and 
$T_0\left(\partial T^\nu(D)\right)\cap J(0)T_0\left(\partial T^{\nu}(D)\right)=\{z_2=0\}$. 
Then by considering the family of vectors $\left(1,0\right)$ at base points $\left(0,t\right)$ for $t\neq 0$ small enough, 
we obtain a family of pseudoholomorphic discs $u_{t}$ such that 
$u_t\left(0\right)=\left(0,t\right)$ and $d_0u_t\left(\partial/\partial_x\right)=\left(0,1\right)$. 
Due to the parameters dependance of the solution to the $J^{\nu}$-holomorphy equation, 
we straighten these discs into the lines $\{z_2=t\}$. Next we  consider a transversal foliation 
by  pseudoholomorphic discs passing through $\left(t,0\right)$ and $\left(t,-\delta_{\nu}\right)$ 
for $t$ small enough and we straighten these lines into $\{z_1=c\}$. This leads to  the 
desired diffeomorphism $\Phi^{\nu}$ of $\R^4$. 

\vspace{0.5cm}

Now, we need to remove harmonic terms from the polynomial
$$\displaystyle \sum_{k=2l_{\nu}}^{2m-1}P_{k}\left(z_1,\overline{z_1},p_{\nu}^*\right).$$
So we consider a biholomorphism (for the standard structure) of $\C^2$ with the following form: 
$$\varphi_{\nu}\left(z_1,z_2\right):=\left(z_1,z_2+\displaystyle \sum_{k=2l_{\nu}}^{2m-1}\Re e \left(c_{k,\nu}z_1^k\right)\right),$$
where $c_{k,\nu}$ are well chosen complex numbers. Then the diffeomorphism $\Phi_\nu:=\varphi_{\nu}\circ\Phi^{\nu}$ satisfies: 
\begin{enumerate}
\item $\Phi_{\nu}\left(p_{\nu}^{*}\right)=0$ and $\Phi_{\nu}\left(p_{\nu}\right)=\left(0,-\delta_{\nu}\right)$. 
\item $\Phi_{\nu}$ converges to $Id~:\R^{4} \rightarrow \R^{4}$ on any compact subset of $\R^{4}$ in the $\mathcal{C}^{2}$ 
sense.
\item  If we denote by $D_{\nu}:=\Phi_{\nu}\left(D\cap U\right)$ the domain with the defining function  
$\rho_{\nu}:=\rho\circ \left(\Phi_{\nu}\right)^{-1}$, then $\rho_{\nu}$ is given by: 
$$\rho_{\nu}\left(z_1,z_2\right)= \Re ez_2+\displaystyle \sum_{k=2l_{\nu}}^{2m-1}P_{k}^{*}\left(z_1,\overline{z_1},p_{\nu}^*\right)+
P_{2m}\left(z_1,\overline{z_1},p_{\nu}^*\right)+O\left(|z_1|^{2m+1}+|z_2|\|z\|\right),$$
where the polynomial 
$$\displaystyle \sum_{k=2l_{\nu}}^{2m-1}P_{k}^{*}\left(z_1,\overline{z_1},p_{\nu}^*\right)$$ 
does not contain any harmonic terms. Moreover the polynomial $P_{2l_{\nu}}^*$  is not idencally zero. 
Moreover, generically, $J_{\nu}:=\left(\Phi_{\nu}\right)_*J$ is no more diagonal.
\end{enumerate}

\vspace{0.5cm}

Since the origin is a boundary point of D'Angelo type $2m$ for $D$, it follows that, denoting by $P_{2m}^*$ the 
nonharmonic part of $P_{2m}$, we have $P_{2m}^*\left(.,0\right)=H_{2m}^*\neq 0$, where $H_{2m}^*$ is the
nonharmonic part of $H_{2m}$. This allows to define for large $\nu$: 
\begin{equation}\label{eqtau}
\tau\left(p_{\nu}^*,\delta_{\nu}\right):= 
\displaystyle \min_{k=2l_\nu,\cdots,2m}\left(\frac{\delta_{\nu}}{\|P_k^*\left(.,p_{\nu}^*\right)\|}\right)^{\frac{1}{k}}.
\end{equation} 
Moreover  the following inequalities hold: 
\begin{equation}\label{ineqq}
\frac{1}{C}\delta_{\nu}^{\frac{1}{2}}\leq \tau\left(p_{\nu}^*,\delta_{\nu}\right)\leq C\delta_{\nu}^{\frac{1}{2m}},
\end{equation}
where $C$ is a positive constant. The right inequality comes from the fact that 
$\|P_{2m}^*\left(.,p_{\nu}^*\right)\|\geq C_1>0$ for large $\nu$. And the left one comes the fact that there exists a positive
constant $C_2$ such that for every $2l_\nu \leq k \leq  2m$ we have $\|P_k^*\left(.,p_{\nu}^*\right)\|\leq C_2$.  

\vspace{0.5cm}

Now we consider the nonisotropic dilation $\Lambda_{\nu}$ of $\C^2$: 
$$\Lambda_{\nu}~: \left(z_1,z_2\right)
\mapsto \left(\tau\left(p_{\nu}^*,\delta_{\nu}\right)^{-1}z_1,\delta_{\nu}^{-1} z_2\right).$$ 

We set $\displaystyle \tilde{D_{\nu}}:=\Lambda_{\nu}\left(D_{\nu}\right)$ the domain admitting the defining function  
$\displaystyle \tilde{\rho_\nu}:=\delta_{\nu}^{-1}\rho_{\nu}\circ \Lambda_{\nu}^{-1}$ and   
$\displaystyle \tilde{J_{\nu}}:= \left(\Lambda_{\nu}\right)_*\left(J_{\nu}\right)$ the direct image of $J_{\nu}$ under 
 $\Lambda_{\nu}$. 

The next lemma is devoted to describe $\displaystyle (\tilde{D_{\nu}},\tilde{J_{\nu}})$ when passing at the limit.
\begin{lem}\label{lemmod}\mbox{ }

\begin{enumerate}
\item The domain $\tilde{D^{\nu}}$ converges in the sense of local Hausdorff set convergence to 
a (standard) pseudoconvex domain $\tilde{D}=\{\tilde{\rho}<0\}$, with 
$$\tilde{\rho}\left(z\right) = \Re e  z_2 + P\left(z_1,\overline{z_1}\right),$$
where $P$ is a nonzero subharmonic polynomial of degree smaller than or equal to $ 2m$ which admits a nonharmonic part. 

\item  \textbf {In case the origin is of D'Angelo type four for $D$}, the sequence of almost complex structures 
$\tilde{J_{\nu}}$ converges on any compact subsets of $\C^{2}$ in the $\mathcal{C}^{2}$ sense to $J_{st}$.

\end{enumerate}
\end{lem}
\begin{proof}[Proof.]

We first prove part 1.
Due to inequalities (\ref{ineqq}), the defining function of $\tilde{D_{\nu}}$ satisfies:
$$\tilde{\rho_{\nu}}=
\Re ez_2+\displaystyle \sum_{k=2l_{\nu}}^{2m}\delta_{\nu}^{-1}\tau\left(p_{\nu}^*,\delta_{\nu}\right)^kP_{k}^*\left(z_1,\overline{z_1},p_{\nu}^*\right)
+\delta_{\nu}^{-1}\tau\left(p_{\nu}^*,\delta_{\nu}\right)^{2m}P_{2m}\left(z_1,\overline{z_1},p_{\nu}^*\right)+O\left(\tau\left(\delta_{\nu}\right)\right).$$
Passing to a subsequence, we may assume that the polynomial 
$$\displaystyle \sum_{k=2l_{\nu}}^{2m}\delta_{\nu}^{-1}\tau\left(p_{\nu}^*,\delta_{\nu}\right)^kP_{k}^*\left(z_1,\overline{z_1},p_{\nu}^*\right)
+\delta_{\nu}^{-1}\tau\left(p_{\nu}^*,\delta_{\nu}\right)^{2m}P_{2m}\left(z_1,\overline{z_1},p_{\nu}^*\right)$$ 
converges uniformly 
on compact subsets of $\C^{2}$ to a nonzero polynomial $P$ of degree $\leq 2m$ admitting a nonharmonic part.  
Since the pseudoconvexity is invariant under diffeomorphisms, it follows that  the domains
$\tilde{D^{\nu}}$ are  $\tilde{J_{\nu}}$-pseudoconvex, and then passing to the limit, the domain 
$\tilde{D}$ is $J_{st}$-pseudoconvex. Thus the polynomial $P$ is subharmonic.

\vskip 0.7cm  

We next prove part 2.  The  complexification of the almost complex structure $J_{\nu}$ is given  by  
\begin{eqnarray*}
\left(J_{\nu}\right)_{\C}& = & \displaystyle \sum_{l=1}^2\Big(A_{l,l}\left(z\right)dz_l\otimes\frac{\partial}{\partial z_l}+
B_{l,l}\left(z\right)dz_l\otimes\frac{\partial}{\partial \overline{z_l}}+
\overline{B_{l,l}}\left(z\right)d\overline{z_l}\otimes\frac{\partial}{\partial z_l}+\\
&&\\
& & \overline{A_{l,l}}\left(z\right)d\overline{z_l}\otimes\frac{\partial}{\partial \overline{z_l}}\Big)
+ A_{1,2}\left(z\right)dz_1\otimes\frac{\partial}{\partial z_2}+
B_{1,2}\left(z\right)dz_1\otimes\frac{\partial}{\partial \overline{z_2}}+
\\
&&\\
& & \overline{B_{1,2}}\left(z\right)d\overline{z_1}\otimes\frac{\partial}{\partial z_2}+
\overline{A_{1,2}}\left(z\right)d\overline{z_1}\otimes\frac{\partial}{\partial \overline{z_2}},
\end{eqnarray*}
where 
$$\left\{
\begin{array}{lll}
A_{l,l}\left(z\right)&=&i+O\left(\left|z_2+\displaystyle \sum_{k=2}^{3} c_{k,\nu}z_1^k\right|^2\right)  \mbox { for } l=1,2,\\
&&\\
B_{l,l}\left(z\right)&=&O\left(\left|z_2+\displaystyle \sum_{k=2}^{3} c_{k,\nu}z_1^k\right|\right) \mbox { for } l=1,2,\\
&&\\
A_{1,2}\left(z\right)&=&\displaystyle \sum_{k=2}^{3} kc_{k,\nu}z_1^{k-1}O\left(
\left|z_2+\displaystyle \sum_{k=2}^{3} c_{k,\nu}z_1^k\right|^2\right),\\
&&\\
B_{1,2}\left(z\right)&=&\displaystyle \sum_{k=2}^{3} k\left(c_{k,\nu}z_1^{k-1}-\overline{c_{k,\nu}z_1^{k-1}}\right)
O\left(\left|z_2+\displaystyle \sum_{k=2}^{3} c_{k,\nu}z_1^k\right|\right).
\end{array}\right.
$$
By a direct computation,  the complexification of 
$\tilde{J_{\nu}}$ is equal to:
\begin{eqnarray*}
\left(\tilde{J_{\nu}}\right)_{\C}& =& 
\displaystyle \sum_{l=1}^2(A_{l,l}(\Lambda_{\nu}^{-1}(z))dz_l\otimes\frac{\partial}{\partial z_l}+
B_{l,l}(\Lambda_{\nu}^{-1}(z))dz_l\otimes\frac{\partial}{\partial \overline{z_l}}+
\\
&&\\
& &  \overline{B_{l,l}}(\Lambda_{\nu}^{-1}(z))d\overline{z_l}\otimes\frac{\partial}{\partial z_l}+
\overline{A_{l,l}}(\Lambda_{\nu}^{-1}(z))d\overline{z_l}\otimes\frac{\partial}{\partial \overline{z_l}})+ \\
&&\\
& & \tau(p_{\nu}^*,\delta_{\nu})\delta_{\nu}^{-1}A_{1,2}(\Lambda_{\nu}^{-1}(z))dz_1\otimes\frac{\partial}{\partial z_2}+
\tau(p_{\nu}^*,\delta_{\nu})\delta_{\nu}^{-1}B_{1,2}(\Lambda_{\nu}^{-1}(z))dz_1\otimes\frac{\partial}{\partial \overline{z_2}}+\\
&&\\
& & \tau(p_{\nu}^*,\delta_{\nu})\delta_{\nu}^{-1}\overline{B_{1,2}}(\Lambda_{\nu}^{-1}(z))
d\overline{z_1}\otimes\frac{\partial}{\partial z_2}+
\tau(p_{\nu}^*,\delta_{\nu})\delta_{\nu}^{-1}\overline{A_{1,2}}(\Lambda_{\nu}^{-1}(z))
d\overline{z_1}\otimes\frac{\partial}{\partial \overline{z_2}}.
\end{eqnarray*}
According to  (\ref{ineqq}) and since $c_{k,\nu}$ converges to zero when $\nu$ tends to $+\infty$ for $k=2,3$, 
it follows that $\tilde{J_{\nu}}$ converges to $J_{st}$. This proves part $\left(2\right)$.
\end{proof}

\subsection{Complete hyperbolicity in D'Angelo type four condition}

In this subsection we prove Theorem \ref{theo}. 
Keeping notations of the previous subsection; we start by establishing the following lemma which gives a  precise 
localization of pseudoholomorphic discs in boxes.  
\begin{lem}\label{LEM} Assume the origin $\in \partial D$ is a point of D'Angelo type four. 
There are positive constants  $C_0$, $\delta_0$ and $r_0$ such that for any
$0 < \delta < \delta_0$, for any large  $\nu$ and for any
${J}_\nu$-holomorphic disc $g_\nu: \Delta \rightarrow {D}_\nu$
we have~:
$$
g_\nu\left(0\right)=(0,-\delta_{\nu}) \Rightarrow g_\nu\left(r_0 \Delta\right) \subset Q\left(0,C_0 \delta_{\nu}\right)
,
$$
where $Q\left(0,\delta_{\nu}\right):=\{z \in \mathbb C^2: |z_1|\leq  \tau\left(p_{\nu}^*,\delta_{\nu}\right),
|z_2| \leq \delta_{\nu}\}$.
\end{lem}
\begin{proof}{\it Proof of Lemma~\ref{LEM}}. Assume by contradiction that there are  a sequence $(C_\nu)_\nu$ 
that tends to $+\infty$ as $\zeta_\nu$ converges to $0$ in $\Delta$, and ${J}_\nu$-holomorphic discs
$g_\nu: \Delta \rightarrow {D}_\nu$ such that
$ g_\nu\left(0\right)=(0,-\delta_{\nu})$ 
and
$g_\nu\left(\zeta_\nu\right) \not\in Q\left(0,C_\nu \delta_\nu\right)$. 
We consider the nonisotropic dilations of $\C^2$:
$$\Lambda_{\nu}^r~: \left(z_1,z_2\right)
\mapsto \left(r^{\frac{1}{4}}\tau\left(p_{\nu}^*,\delta_{\nu}\right)^{-1}z_1,r\delta_{\nu}^{-1} z_2\right),$$ 
where $r$  is a positive constant to be fixed. 
We set  $h_\nu:= \Lambda_{\nu}^r\circ g_\nu$, $\tilde{\rho_\nu^r}:=r\delta_{\nu}^{-1}\rho_\nu \circ \left(\Lambda_{\nu}^r\right)^{-1}$ 
and  $\tilde{J_\nu^r}:=\left(\Lambda_{\nu}^r\right)_*J_\nu$.
It follows from Lemma~\ref{lemmod} that
$\tilde{\rho_\nu^r}$ converges to 
$$\tilde{\rho}=Re\left(z_2\right) + P\left(z_1,\overline{z_1}\right)$$
uniformly on any compact subset of $\mathbb C^2$ and $\tilde{J_\nu^r}$ converges
to $J_{st}$, uniformly on any compact subset of $\mathbb C^2$. 
 According to the stability of the Kobayashi pseudometric stated in Proposition \ref{propstab}, 
there exist a positive constant $C$ and a neighborhood $V$ of the origin 
 in $\R^4$, such that for every large $\nu$, for  every $q \in \tilde{D_{\nu}} \cap V$ and every $v\in T_q \R^4$:  
$$K_{\left(\tilde{D_{\nu}},\tilde{J_\nu}\right)}\left(q,v\right)\geq C \|v\|.$$ 
Therefore, there exists a constant $C' > 0$ such that 
$$\parallel dh_\nu\left(\zeta\right) \parallel \leq C'$$ 
for any $\zeta \in \left(1/2\right)\Delta$ satisfying $h_\nu\left(\zeta\right) \in \tilde{D_\nu}\cap V'$, with $ V' \subset V$. 
Now we choose the constant $r$ such that  $h_\nu\left(0\right)=\left(0,-r\right) \in Int\left(V'\right)$.
On the other hand, the sequence $\vert h_\nu\left(\zeta_\nu\right) \vert$ tends to $+
\infty$. Denote by $[0,\zeta_\nu]$ the segment 
(in $\mathbb C$) joining the origin and $\zeta_\nu$ and let 
$\zeta_\nu'=r_\nu e^{i\theta_\nu} \in [0,\zeta_\nu]$ be the point closest to the origin such
that  $h_\nu\left([0,\zeta_\nu']\right) \subset \tilde{D_\nu} \cap V$ and 
$h_\nu\left(\zeta_\nu'\right) \in \partial V$. Since $h_\nu\left(0\right)
\in Int\left(V'\right)$, we have  
$$\| h_\nu\left(0\right) - h_\nu\left(\zeta_\nu'\right) \| \geq C''$$
for some constant $C'' > 0$. It follows that: 
$$
\| h_\nu\left(0\right) - h_\nu\left(\zeta_\nu'\right)\| \leq \int_{0}^{r_\nu}
\left\| dh_\nu\left(te^{i\theta_\nu}\right) \right\| dt \leq C'r_\nu\longrightarrow 0.
$$
This contradiction proves Lemma~\ref{LEM}.
\end{proof}

Now we go on the proof of Theorem \ref{theo}.
\begin{proof}[Proof of Theorem \ref{theo}]

Due to the localization of the Kobayashi pseudometric established in   Proposition \ref{thm}, 
it suffices to prove Theorem \ref{theo} in a 
neighborhood $U$ of $q  \in \partial D$. 
Choosing  local coordinates  $z~: U\rightarrow \B \subset \R^{4}$ centered at $q$, 
we may assume that $D\cap U=\{\rho<0\}$ is a $J$-pseudonconvex region of $(\R^{4},J)$, 
that $q=0 \in \partial D$ and  that $J$ satisfies (\ref{eqstr}) and (\ref{eqstrO}).  
We also suppose that the complex tangent space $T_{0}\partial D\cap J(0)T_{0}\partial D$ 
at $0$ of $\partial D$ is given by $\{z_{2}=0\}$. Moreover the defining function $\rho$ is expressed by:    
$$\rho\left(z\right) =  \Re ez_2+H_{2m}\left(z_1,\overline{z_1}\right)+O\left(|z_1|^{2m+1}+|z_2|\|z\|\right)$$

\vspace{0.5cm}

For $p\in D\cap U$ be sufficiently close to the boundary $\partial D$, there exists a unique point 
$p^* \in \partial D \cap U$ such that 
$$p^*=p+(0,\delta),$$
with $\delta>0$. We define an infinitesimal pseudometric $N$ on $D\cap U \subseteq \R^4$ by: 
\begin{equation}
N\left(p,v\right):=\frac{|\left(d_p\Phi_{p^*}v\right)_1|}{\tau\left(p^*,|\rho\left(p\right)|\right)}+
\frac{|\left(d_p\Phi_{p^*}v\right)_{2}|}{|\rho\left(p\right)|},
\end{equation}
for every $p \in D \cap U$ and every $v \in T_{p}\R^{4}$, where  $\Phi_{p^*}$ is defined as 
diffeomorphisms $\Phi^{\nu}$ (of previous subsection) for $p^*$ instead of $p_\nu^*$.

\vspace{0.5cm}

To prove estimate (\ref{est1}) of Theorem \ref{theo}, it suffices to find a positive 
constant $C$ such that for any $J$-holomorphic disc 
$u~: \Delta \rightarrow D\cap U$, we have:  
\begin{equation}\label{eqeq}
N\left(u\left(0\right),d_{0}u\left(\partial/\partial_{x}\right)\right)\leq C.
\end{equation}
Indeed, for a $J$-holomorphic disc $u$ such that $u\left(0\right)=p$ and  $d_{0}u\left(\partial/\partial_{x}\right)=r v$, 
(\ref{eqeq}) leads to  
$$\frac{1}{r} = \frac{N\left(p,v\right)}{N\left(u\left(0\right),
d_{0}u\left(\partial/\partial_{x}\right)\right)}\geq \frac{N\left(p,v\right)}{C}.$$

\vskip 0.5cm

Suppose by contradiction that (\ref{eqeq}) is not true, that is, there is a sequence of $J$-holomorphic discs 
$u_\nu~: \Delta \rightarrow D\cap U$ such that $N\left(u_\nu\left(0\right),
d_{0}u_\nu\left(\partial/\partial_{x}\right)\right)\geq \nu^{2}.$
Then we consider  a sequence $\left(y_{\nu}\right)_\nu$ of points in  $\overline{\Delta_{1/2}}$ such that:
\begin{enumerate}
\item $\displaystyle |y_{\nu}|\leq \frac{2\nu}{N\left(u_\nu\left(y_{\nu}\right),d_{y_{\nu}}u_\nu\left(\partx\right)\right)}$ , \\
\item  $\displaystyle N\left(u_\nu\left(y_{\nu}\right),d_{y_{\nu}}u_\nu\left(\partx\right)\right)\geq \nu^{2}$, and \\
\item $\displaystyle y_{\nu}+\Delta_{\nu/N\left(u_\nu\left(y_{\nu}\right),d_{y_{\nu}}u_\nu\left(\partx\right)\right)}\subseteq 
 \overline{\Delta_{1/2}}$ for sufficiently large $\nu$.\\  
\end{enumerate}
This allows to define a sequence of  $J$-holomorphic discs $g_{\nu}~: \Delta_{\nu}\rightarrow D\cap U$  by 
$$g_{\nu}\left(\zeta\right):=u_\nu\left(y_{\nu}+\frac{\zeta}{2N\left(u_\nu\left(y_{\nu}\right),
d_{y_{\nu}}u_\nu\left(\partx\right)\right)}\right).$$
Consider the sequence $g_\nu=u_\nu(y_\nu)$ in $D\cap U$. Since $|y_\nu|\leq 2/\nu$ and since the $\mathcal{C}^1$ norm of 
any $J$-holomorphic disc $u_\nu$ is uniformally bounded it follows that $g_\nu(0)$ converges to the origin.

\vspace{0.5cm}

We apply the scaling method to the sequence $g_{\nu}\left(0\right)$. We denote by 
$g_{\nu}\left(0\right)^*$ the boundary point given by  
$g_{\nu}\left(0\right)^*:=g_{\nu}\left(0\right)+\left(0,\delta_{\nu}\right)$. We set the scaled 
disc $\tilde{g_{\nu}}:=\Lambda_{\nu}\circ \Phi_{\nu}\circ g_{\nu}$, where diffeomorphisms 
$\Lambda_{\nu}$ and $\Phi_{\nu}$ are define in the subsection about the scaling method. 
In order to extract from $\tilde{g_{\nu}}$ a subsequence which converges to a 
Brody curve, we need the following Lemma.
\begin{lem}\label{lemlee}  There is a positive constant $r_0$ such that:
\begin{enumerate} 
\item There exists a positive constant $C_1$  such that 
\begin{equation}\label{eqincinc}
\tilde{g_{\nu}}\left(r_0\Delta_\nu\right)\subset  \Delta_{C_1}\times\Delta_{C_1}.
\end{equation}
\item There exists a positive constant $C_2$ such that for every  large $\nu$ we have~: 
\begin{equation}\label{estder}
\|d\tilde{g_{\nu}}\|_{\mathcal{C}^{0}\left(r_0\Delta_{\nu}\right)}\leq C_2.
\end{equation}
\end{enumerate}
\end{lem}
\begin{proof}[Proof.]
We prove the first part. We define a $J_{\nu}$-holomorphic disc 
$h_{\nu}\left(\zeta\right):=\Phi_\nu\circ g_{\nu}\left(\nu\zeta\right)$ from the unit disc 
$\Delta$ to $D_{\nu}$. According to Lemma \ref{LEM}, 
since $h_{\nu}\left(0\right)=\Phi_\nu\circ g_{\nu}\left(0\right)= (0,-\delta_\nu)$, we have 
$$h_{\nu}\left(r_0\Delta\right)\subseteq Q\left(0,C_0\delta_\nu\right)$$ 
for some positive constants $r_0$ and $C_0$. Hence 
$$\Phi_\nu \circ g_{\nu}\left(r_0\Delta_\nu\right)\subseteq Q\left(0,C_0\delta_\nu\right).$$ 
After dilations, this leads to (\ref{eqincinc}).

\vspace{0.7cm}

Then we prove the second part. According to Lemma \ref{lemmod}, the sequence of almost complex structures 
 $\tilde{J_{\nu}}$ converges on any compact subsets of $\C^{2}$ 
in the $\mathcal{C}^{2}$ sense to $J_{st}$. 
Then for sufficiently large $\nu$, the norm 
$\|\tilde{J_{\nu}}-J_{st}\|_{\mathcal{C}^{1}(\overline{\Delta_{C_1}\times\Delta_{C_1}})}$ is as small as necessary.
So for large  $\nu$, and due to Proposition $2.3.6$ of J.-C.Sikorav in 
\cite{jcsi} there exists  $C_2>0$ such that (\ref{estder}) holds.    
\end{proof}

\vspace{0.5cm}

Hence according to Lemmas \ref{lemmod} and \ref{lemlee} we may extract from $\tilde{g_{\nu}}$ 
a subsequence, still denoted by $\tilde{g_{\nu}}$ wich converges  in  $\mathcal{C}^{1}$ topology
to a standard complex line 
$$\tilde{g}~: \C \rightarrow \left(\{Re z_2 + P\left(z_1,\overline{z_1}\right)<0\},J_{st}\right).$$ 
The polynomial $P$ is subharmonic and contains a nonharmonic part; this implies that  
the domain $\left(\{Re z_2 + P\left(z_1,\overline{z_1}\right)<0\},J_{st}\right)$ is Brody hyperbolic and so  
the complex line $\tilde{g}$ is constant. To obtain a contradiction, 
we prove that the derivative of $\tilde{g}$ at the origin is nonzero:   
$$\frac{1}{2}=N\left(g_\nu\left(0\right),d_{0}g_{\nu}\left(\partial/\partial_{x}\right)\right)=
\frac{|\left(d_{0}\left(\Phi_\nu \circ g_{\nu}\right)\left(\partial/\partial_{x}\right)\right)_1|}{\tau\left(g_{\nu}\left(0\right)^*,|\rho\left(g_{\nu}\left(0\right)\right)|\right)}+
\frac{|\left(d_{0}\left(\Phi_\nu \circ g_{\nu}\right)\left(\partial/\partial_{x}\right)\right)_2|}{|\rho\left(g_{\nu}\left(0\right)\right)|}.$$
Since  $|\rho\left(g_{\nu}\left(0\right)\right)|$ is equivalent to $\delta_\nu$, 
it follows that for some
positive constant $C_3$ and for large $\nu$, we have:  
$$\frac{1}{2}\leq C_3\left(\frac{|\left(d_{0}\left(\Phi_\nu \circ g_{\nu}\right)
\left(\partial/\partial_{x}\right)\right)_1|}{\tau\left(g_{\nu}\left(0\right)^*,\delta_{\nu}\right)}+
\frac{|\left(d_{0}\left(\Phi_\nu \circ g_{\nu}\right)\left(\partial/\partial_{x}\right)\right)_{2}|}{\delta_\nu}\right)
=C_3\|d_{0}\tilde{g_\nu}\left(\partial/\partial_{x}\right)\|_1.$$
Since $\tilde{g_\nu}$ converges to $\tilde{g}$  in the $\mathcal{C}^{1}$ sense, 
it follows that $d_{0}\tilde{g}\left(\partial/\partial_{x}\right)\neq 0$, providing a contradiction. This achieves the proof of 
Theorem \ref{theo}.
\end{proof}

Estimate (\ref{est2}) of the Kobayashi pseudometric allows to study the completness of the Kobayashi pseudodistance $D$. 
\begin{cor}\label{corcomp}
Let $D=\{\rho<0\}$ be a relatively compact domain of finite D'Angelo type less than or equal to four in an almost 
complex manifold $\left(M,J\right)$ of dimension four, where $\rho$ is a defining function of $D$, $J$-plurisubharmonic 
in a neighborhood of $\overline{D}$. Assume that $\left(M,J\right)$ admits a global strictly $J$-plurisubharmonic function. 
Then $\left(D,J\right)$ is complete hyperbolic.
\end{cor}

\begin{proof}[Proof.]
The fact that  $\left(M,J\right)$ admits a global strictly $J$-plurisubharmonic function 
 and estimate (\ref{ee1}) of Proposition \ref{thm} leads to  the Kobayashi hyperbolicity of $D$. Then 
estimate (\ref{est2}) of the Kobayashi pseudometric stated in Theorem \ref{theo} gives the completness of the 
metric space  $\left(D,d_{\left(D,J\right)}\right)$ by a classical integration argument.  
\end{proof}

\subsection{Regions with noncompact automorphisms group}

The next corollary is devoted to regions with noncompact automorphisms group.
\begin{cor}\label{coraut}
Let $D=\{\rho<0\}$ be a relatively compact  domain  in a four dimensional  
almost complex manifold $\left(M,J\right)$ of finite 
D'Angelo type less than or equal to four. 
Assume that $\rho$ is a 
$\mathcal{C}^{2}$ defining function of $D$, $J$-plurisubharmonic on a neighborhood of $\overline{D}$. 
If there is an automorphism of $D$ with orbit accumulating at a boundary point then there exists a polynomial $P$ of degree at 
most four, without harmonic part such that $\left(D,J\right)$ is biholomorphic to 
$\left(\{\Re e z_2+P\left(z_1,\overline{z_1}\right)<0\},J_{st}\right)$.
\end{cor}

If the domain $D$ is a relatively compact  strictly $J$-pseudoconvex domain with noncompact automorphisms group then 
$\left(D,J\right)$ is biholomorphic to a model domain. This was proved  by H.Gaussier and A.Sukhov in \cite {ga-su} 
in  dimension four and by K.H.Lee in \cite{le} in arbitrary (even) dimension.  

\begin{proof}[Sketch of the proof.]

We suppose that for some point $p_0 \in D$, 
there is a sequence $f_\nu$ of automorphisms of $(D,J)$ such that $p_\nu:=f_\nu\left(p_0\right)$ 
converges to $0 \in \partial D$.
We apply the scaling method to the sequence $p_\nu$. Still keeping notations of subsection $4.1$, we set 
$$F_\nu:=\Lambda_\nu \circ \Phi_\nu \circ f_\nu : f_\nu^{-1}\left(D\cap U\right) \rightarrow \tilde{D_\nu}.$$ 
This sequence of biholomorphisms is such that  : 
\begin{enumerate}
\item $\left( f_\nu^{-1}\left(D\cap U\right)\right)_\nu$ converges to $D$.
\item  $\tilde{D_{\nu}}$ converges to a pseudoconvex domain $\tilde{D}=\{ Re z_2 + P\left(z_1,\overline{z_1}\right)<0\}$, 
where $P$ is a nonzero subharmonic polynomial of degree $\leq 4$ which contains a nonharmonic part. 
Changing $\tilde{D}$ by applying a  standard biholomorphism if necessary, 
we may suppose that  $P\left(z_1,\overline{z_1}\right)$ is without harmonic terms.
\item For any compact subset $K \subset D$, the sequence $\left(\|F_\nu\|_{\mathcal{C}^1\left(K\right)}\right)_\nu$ is bounded.

\end{enumerate}
Hence, we may extract from $\left(F_\nu\right)_\nu$ a subsequence converging, 
on any  compact subset of $D$ in the $\mathcal {C^1}$ sense, to a $\left(J,J_{st}\right)$-holomorphic map 
$F : D \longrightarrow \bar{\tilde{D}}$.
Finally $F$ is a $\left(J,J_{st}\right)$-biholomorphism from $D$ to $\tilde{D}$.
\end{proof}

\subsection{Nontangential approach in the general setting}

In this subsection, refering to I.Graham \cite{gra}, we give a sharp estimate of the Kobayashi pseudometric of a 
pseudoconvex region in a cone with vertex at a boundary point of arbitrary finite D'Angelo type. 
We denote by $\Lambda:=\{-\Re e z_2>k\|z\|\}$, where $0<k<1$, the cone with vertex at the origin and axis 
the negative real $z_2$ axis. 
 
\begin{theo}\label{theo3}
Let $D=\{\rho<0\}$ be a   domain  of finite D'Angelo type  in $\left(\R^4,J\right)$, where 
$$\rho\left(z_1,z_2\right)=\Re ez_2+H_{2m}\left(z_1,\overline{z_1}\right)+O\left(|z_1|^{2m+1}+|z_2|\|z\|\right),$$
is a $\mathcal{C}^{2}$ defining function of $D$, $J$-plurisubharmonic on a neighborhood of $\overline{D}$. 
We suppose that $H_{2m}$ is a homogeneous subharmonic polynomial of degree $2m$ admitting a nonharmonic part.  
Then there exists a positive constant $C$ such that for every $p \in D\cap \Lambda$ and every $v \in T_{p}M$: 
\begin{equation*}
K_{\left(D,J\right)}\left(p,v\right)\geq C\left(\frac{|v_1|}{|\rho\left(p\right)|^{\frac{1}{2m}}}+
\frac{|v_{2}|}{|\rho\left(p\right)|}\right).
\end{equation*}
\end{theo} 

Before proving Theorem \ref{theo3} we need the following crucial lemma.
%which implies that when doing the scaling method we only dilate coordinates.  
\begin{lem}\label{lemuti} There exist a neighborhood $U$ of the origin and a positive constant $C$ such that 
if $p \in D \cap U\cap  \Lambda$ then 
$$ p\in  \left\{z \in \C^2~: |z_1|<C_1{\rm dist}\left(p,\partial D\right)^{\frac{1}{2m}}, 
|z_2|<C_1{\rm dist}\left(p,\partial D\right)\right\}.$$
\end{lem}
\begin{proof}[Proof.]
According to the fact that ${\rm dist}\left(z,\partial D\right)$ is equivalent to 
$|\rho\left(z\right)|=-\Re e z_2 +O\left(\|z\|^2\right)$ and to the 
definition of the cone $\Lambda$, we have: 
$$\lim_{z\rightarrow 0, z\in D\cap\Lambda} \frac{-\Re e z_2}{{\rm dist}\left(z,\partial D\right)}=1.$$

This implies the existence of a positive constant $C_1$  such that 
$$\|p\|<-\frac{1}{k}\Re e p_2\leq C_1 {\rm dist}\left(p,\partial D\right),$$ 
whenever  $p\in D\cap \Lambda$ is sufficiently close to the origin.
Thus 
$$  p \in \left\{z \in \C^2~: |z_1|< C_1 {\rm dist}\left(p,\partial D\right)^{\frac{1}{2m}}, 
|z_2|< C_1{\rm dist}\left(p,\partial D\right)\right\},$$
for $p\in D\cap \Lambda$ sufficiently close to the origin.

\end{proof}

The proof of Theorem \ref{theo3} is  similar and easier than  proof of Theorem \ref{theo}. For convenience, we write it. 
\begin{proof}[Proof of Theorem \ref{theo3}.]
Let $U$ be a neighborhood of the origin. We define an infinitesimal pseudometric $N$ on $D\cap U \subseteq \R^4$ by: 
\begin{equation*}
N\left(p,v\right):=\frac{|v_1|}{|\rho\left(p\right)|^{\frac{1}{2m}}}+\frac{|v_2|}{|\rho\left(p\right)|},
\end{equation*}
for every $p \in D \cap U$ and every $v \in T_{p}\C^{2}$.

\vspace{0.5cm}

We have to find a positive constant $C$ such that for every $J$-holomorphic disc 
$u~: \Delta \rightarrow D\cap U$, such that if  $u\left(0\right)\in \Lambda$ then:  
\begin{equation*}
N\left(u\left(0\right),d_{0}u\left(\partial/\partial_{x}\right)\right)\leq C.
\end{equation*}
Suppose by contradiction that this inequality is not true, that is, there exists a sequence of $J$-holomorphic discs 
$u_\nu~: \Delta \rightarrow D\cap U$ such that 
$$u_\nu\left(0\right)\in \Lambda \mbox{ }\mbox{ }\mbox{ and } \mbox{ }\mbox{ } 
N\left(u_\nu\left(0\right),d_{0}u_\nu\left(\partial/\partial_{x}\right)\right)\geq \nu^{2}.$$ 
Then consider  a sequence $\left(y_{\nu}\right)_\nu$ of points in  $\overline{\Delta_{1/2}}$ such that 
\begin{enumerate}
\item $\displaystyle |y_{\nu}|\leq \frac{2\nu}{N\left(u_\nu\left(y_{\nu}\right),d_{y_{\nu}}u_\nu\left(\partx\right)\right)}$,\\ 
\item  $\displaystyle N\left(u_\nu\left(y_{\nu}\right),d_{y_{\nu}}u_\nu\left(\partx\right)\right)\geq \nu^{2}$, and \\
\item $\displaystyle y_{\nu}+\Delta_{\nu/N\left(u_\nu\left(y_{\nu}\right),d_{y_{\nu}}u_\nu\left(\partx\right)\right)}
\subseteq  \overline{\Delta_{1/2}}$ for sufficiently large $\nu$. \\ 
\end{enumerate}
Then we define a sequence of  $J$-holomorphic discs $g_{\nu}~: \Delta_{\nu}\rightarrow D\cap U$  by 
$$g_{\nu}\left(\zeta\right):=u_\nu\left(y_{\nu}+\frac{\zeta}{2N\left(u_\nu\left(y_{\nu}\right),
d_{y_{\nu}}u_\nu\left(\partx\right)\right)}\right).$$

\vspace{0,5cm}

For large $\nu$, we have $g_{\nu}\left(0\right)=u_\nu\left(y_{\nu}\right)$ in $D\cap U\cap \Lambda$ and 
$g_{\nu}\left(0\right)$ converges to the origin. Set 
$$\delta_\nu:={\rm dist}\left(g_\nu\left(0\right),\partial D\right),$$
and consider the following dilations of $\C^2$: 
$$\Lambda_{\nu}~: \left(z_1,z_2\right)\mapsto \left(\delta_{\nu}^{\frac{-1}{2m}}z_1,\delta_{\nu}^{-1} z_2\right).$$ 
In order to extract from $\Lambda_{\nu}\circ g_{\nu}$ a subsequence which converges to a 
Brody curve, we need the following Lemma.
\begin{lem}\label{lemlee2}  There exists a positive constant $r_0$ such that:
\begin{enumerate} 
\item there exists a positive constant $C_1$ 
such that:
\begin{equation}\label{lemlee4}
\Lambda_{\nu}\circ g_{\nu}\left(r_0\Delta_\nu\right)\subset  \Delta_{C_1}\times\Delta_{C_1},
\end{equation}
\item there is a positive constant $C_2$ such that for every large  $\nu$ we have~: 
\begin{equation}\label{lemlee3}
\|d\left(\Lambda_{\nu}\circ g_{\nu}\right)\|_{\mathcal{C}^{0}\left(r_0\Delta_{\nu}\right)}\leq C_2.
\end{equation}
\end{enumerate}
\end{lem}
\begin{proof}[Proof.]
We first prove (\ref{lemlee4}). We define a new $J$-holomorphic disc 
$h_{\nu}\left(\zeta\right):=g_{\nu}\left(\nu\zeta\right)$ from the unit disc 
$\Delta$ to $D_{\nu}$. According to Lemma \ref{lemuti}, we have 
$$h_{\nu}\left(0\right)=g_{\nu}\left(0\right) \in \{z \in \C^2:
|z_1|\leq C_1\delta_{\nu}^{\frac{1}{2m}}, |z_2|\leq C_1\delta_{\nu}\}.$$
This implies:  
$$h_{\nu}\left(r_0\Delta\right)\subseteq \{z \in \C^2~: |z_1|\leq C_0\delta_{\nu}^{\frac{1}{2m}},
|z_2| < C_0\delta_{\nu}\},$$
 for positive constants $r_0$ and $C_0$, since  Lemma \ref{LEM} is true if we replace $\tau\left(p_{\nu}^*,\delta_{\nu}\right)$
by $\delta_{\nu}^{\frac{1}{2m}}$.
Hence 
$$g_{\nu}\left(r_0\Delta_\nu\right)\subseteq \{z \in \C^2~: |z_1|< C_0\delta_{\nu}^{\frac{1}{2m}},
|z_2| \leq C_0\delta_{\nu}\}.$$ 
After dilations, this leads to  (\ref{lemlee4}).

\vspace{0.7cm}

The proof of (\ref{lemlee3}) is similar to (\ref{estder}) of Lemma \ref{lemlee}, since the sequence of structures $\left(\Lambda_{\nu}\right)_*J$ 
converges on any compact subset of $\C^{2}$
in the $\mathcal{C}^{1}$ sense to $J_{st}$ because $J$ is diagonal.  
\end{proof}

Hence according to Lemma \ref{lemlee2} we may extract from $\Lambda_\nu\circ g_\nu$ 
a subsequence, still denoted by $\Lambda_{\nu}\circ g_{\nu}$ wich converges in the $\mathcal{C}^{1}$ sense
to a standard complex line $\tilde{g}~: \C \rightarrow \left(\{Re z_2 + H_{2m}\left(z_1,\overline{z_1}\right)<0\},J_{st}\right)$, 
where the domain $\left(\{Re z_2 + P\left(z_1,\overline{z_1}\right)<0\},J_{st}\right)$ 
is Bro\-dy hyperbolic since  $H_{2m}\left(z_1,\overline{z_1}\right)$ 
contains a nonharmonic part. Then the standard complex line $\tilde{g}$ is constant. To obtain a contradiction, 
we prove that the derivative of $\tilde{g}$ is nonzero:    
$$\frac{1}{2}=N(g_\nu(0),d_{0}g_{\nu}(\partial/\partial_{x}))=
\frac{|(d_{0}g_{\nu}(\partial/\partial_{x}))_1|}{|\rho(g_{\nu}(0))|^{\frac{1}{2m}}}+
\frac{|(d_{0}g_{\nu}(\partial/\partial_{x}))_2|}{|\rho(g_{\nu}(0))|}.$$
Since  $|\rho\left(g_{\nu}\left(0\right)\right)|$ is equivalent to $\delta_\nu$, 
it follows that for some
positive constant $C_3$ we have for large $\nu$:  
$$\frac{1}{2}\leq C_3\left(\frac{|(d_{0}(g_{\nu})(\partial/\partial_{x}))_1|}
{\delta_{\nu}^\frac{1}{2m}}+
\frac{|(d_{0}(g_{\nu})(\partial/\partial_{x}))_{2}|}{\delta_\nu}\right)
=C_3\|d_{0}(\Lambda_\nu \circ g_\nu)(\partial/\partial_{x})\|_1.$$
This provide a contradiction.
\end{proof}

\section{Appendix : Convergence of the structures involved by the scaling method.} 
%\bigskip \noindent {\bf 1. About the assumption of D'Angelo type 4.} 
\bigskip 
In this appendix, we prove that, generically, the convergence of the sequence of structures involved by the scaling method  
to the standard structure  $J_{st}$ occurs only on a neighborhood of  boundary points of D'Angelo type less than or equal to four. 

\vspace{0.2cm}

Let  $D=\{\rho<0\}$ be  a pseudoconvex region of finite D'Angelo type $2m$ in $\R^4$, where $\rho$ has the following 
expression on a neighborhood $U$ of the origin:   
$$\rho\left(z_1,z_2\right)=\Re ez_2+H_{2m}\left(z_1,\overline{z_1}\right)+O\left(|z_1|^{2m+1}+|z_2|\|z\|\right),$$
where $H_{2m}$ is a homogeneous subharmonic polynomial of degree $2m$ admitting a nonharmonic part. 
Assume that $p_{\nu}$ is a sequence of points in $D\cap U$ converging to the origin,  
and, for large $\nu$, consider the sequence of diffeomorphisms $\Phi^{\nu}: \R^4\rightarrow \R^4$ given in the scaling method. 
We suppose that the function  $\rho^{\nu}=\rho\circ \left(\Phi^{\nu}\right)^{-1}$ is given by: 
$$\rho^{\nu}\left(z_1,z_2\right)= \Re e z_2+\Re e \left(\alpha_\nu z_1^2\right)+\beta\nu|z_1|^2
+\displaystyle \sum_{k=3}^{2m}P_{k}\left(z_1,\overline{z_1},p_{\nu}^*\right)+O\left(|z_1|^{2m+1}+|z_2|\|z\|\right).$$
Moreover the structure $J^\nu:=\left(\Phi^{\nu}\right)_*J$ satisfies (\ref{eqstr}) and (\ref{eqstrO}).  
To fix notations, we set: 
\begin{equation*}
J^{\nu}=\left(\begin{array}{ccccc} 
a_1^\nu& b_1^\nu & 0 & 0 \\
c_1^\nu & -a_1^\nu & 0 & 0\\
0  & 0 & a_2^\nu & b_2^\nu\\
0 & 0 & c_2^\nu & -a_2^\nu\\ 
\end{array}\right).
\end{equation*}

\vspace{0,5cm}

Now, consider the following diffeomorphism of $\R^4$ defined by:
\begin{equation}\label{eqdif}
\Psi_\nu^{-1}\left(x_1,y_1,x_2,y_2\right)=\left(x_1+R_{1,\nu},y_1+S_{1,\nu},x_2+R_{2,\nu},y_2+S_{2,\nu}\right)
\end{equation}
converging to the identity and such that $d_0\Psi_\nu^{-1}=Id$. We suppose that $R_{k,\nu}$ and $S_{k,\nu}$, for $k=1,2$ 
are  real functions depending smoothly on $x_1,y_1$ and $y_2$ and that $R_{2,\nu}$ and $S_{2,\nu}$ are given by: 
\begin{equation}\label{eqdif1}
\left\{
\begin{array}{lll}  
R_{2,\nu}& = & -\alpha_\nu x_1^2+ \alpha_\nu y_1^2+O\left(|z_1|^3+y_2^2+|y_2|\|z\|\right), \\
& & \\
S_{2,\nu}& = & -2 \alpha_\nu x_1y_1+O\left(|z_1|^3+y_2^2+|y_2|\|z\|\right).\\
\end{array}
\right.
\end{equation}
We write:
\begin{equation}\label{eqdif2}
\left\{
\begin{array}{lll}  
R_{1,\nu}& = & r_{5,\nu}x_1^2+r_{6,\nu}x_1y_1+r_{7,\nu}y_1^2+r_{1,\nu}x_1^3+r_{2,\nu}x_1^2y_1+r_{3,\nu}x_1y_1^2+\\
& & \\
& & r_{4,\nu}y_1^3+O\left(|z_1|^4+y_2^2+|y_2|\|z\|\right)\\
& & \\
S_{1,\nu}& = & s_{5,\nu}x_1^2+s_{6,\nu}x_1y_1+s_{7,\nu}y_1^2+s_{1,\nu}x_1^3+s_{2,\nu}x_1^2y_1+s_{3,\nu}x_1y_1^2+\\
& & \\
& & s_{4,\nu}y_1^3+O\left(|z_1|^4+y_2^2+|y_2|\|z\|\right).\\
\end{array}
\right.
\end{equation}
%\begin{eqnarray*} 
%R_{1,\nu}& = & r_{5,\nu}x_1^2+r_{6,\nu}x_1y_1+r_{7,\nu}y_1^2+O\left(|z_1|^3+y_2^2+|y_2|\|z\|\right), \mbox{ and, }\\
%S_{1,\nu}& = & s_{5,\nu}x_1^2+s_{6,\nu}x_1y_1+s_{7,\nu}y_1^2+O\left(|z_1|^3+y_2^2+|y_2|\|z\|\right).\\
%\end{eqnarray*} 
It follows that:
\begin{eqnarray*}
\rho_\nu\circ \Psi_\nu^{-1}\left(z_1,z_2\right)& =& \Re e z_2+\beta_\nu|z_1^2|+ 
\displaystyle \sum_{k=3}^{2m}P_{k}'\left(z_1,\overline{z_1},\nu\right)+O\left(|z_1|^{2m+1}+|z_2|\|z\|\right).\\
\end{eqnarray*}

\vspace{0,5cm}
Then we define
$$\tau_\nu:=\min\left(\left(\frac{\delta_\nu}{|\beta_\nu|}\right)^{\frac{1}{2}},
\min_{k=3,\cdots,2m-1}\left(\frac{\delta_\nu}{\|P_k'\left(.,\nu\right)\|}\right)^{\frac{1}{k}},\delta_\nu^{\frac{1}{2m}}\right).$$
And we consider the following anisotropic dilations of $\C^2$:
$$\Lambda_\nu\left(z_1,z_2\right):=\left(\tau_\nu^{-1}z_1,\delta_\nu^{-1}z_2\right).$$
If we write $J_{\nu}:=\left(\Psi_{\nu}\right)_*J^\nu$ as:
\begin{equation*}
J_{\nu}=\left(\begin{array}{ccccc} 
J_{1,\nu}& B_{1,\nu}\\
C_{1,\nu} & J_{2,\nu}\\
\end{array}\right) \mbox { with }  \mbox { } 
C_{1,\nu}:=\left(\begin{array}{ccccc} 
\left(J_\nu\right)_1^3& \left(J_\nu\right)_2^3\\
\left(J_\nu\right)_1^4& \left(J_\nu\right)_2^4\\
\end{array}\right),
\end{equation*}
then we have: 
\begin{equation*}
\left(\Lambda_{\nu}\right)_*J_{\nu}\left(z\right)=\left(\begin{array}{ccccc} 
J_{1,\nu}\left(\tau_\nu z_1,\delta_\nu z_2\right)& \tau_\nu^{-1} \delta_\nu B_{1,\nu}\left(\tau_\nu z_1,\delta_\nu z_2\right)\\
\tau_\nu \delta_\nu^{-1} C_{1,\nu}\left(\tau_\nu z_1,\delta_\nu z_2\right) & J_{2,\nu}\left(\tau_\nu z_1,\delta_\nu z_2\right)\\
\end{array}\right).
\end{equation*}

\vspace{0,5cm}

We have generically the following situation:
\begin{prop}
The sequence of structures $\left(\Lambda_{\nu}\right)_*J_{\nu}$ converges to the standard structure $J_{st}$ 
if and only if the D'Angelo type of the origin is less than or equal to four.
\end{prop}
\begin{proof}[Proof.]
We notice that $\left(\Lambda_{\nu}\right)_*J_{\nu}$ converges to $J_{st}$ if and only if  
$C_{1,\nu}=O\left(|z_1|^{2m-1}\right)+O\left(|z_2|\right)$. 
Indeed  if  $C_{1,\nu}=O\left(|z_1|^{2m-1}\right)+O\left(|z_2|\right)$ then 
$$\tau_\nu \delta_\nu^{-1} C_{1,\nu}\left(\tau_\nu z_1,\delta_\nu z_2\right)=\tau_\nu^{2m}\delta_\nu^{-1}O|z_1|^{2m}+
\tau_\nu^{2m}O|z_1|^{2m},$$
which converges to  the zero $2$ by $2$ matrix since $\tau_\nu \leq \delta_\nu^{\frac{1}{2m}}$ and 
since $C_{1,\nu}$ tends to the zero $2$ by $2$ matrix. 
Conversely  if $C_{1,\nu}=O\left(|z_1|^{k}\right)+O\left(|z_2|\right)$, with $k<2m-1$, 
then  $\left(\Lambda_{\nu}\right)_*J_{\nu}$ converges to a polynomial 
integrable structure $\tilde{J}=J_{st}+O|z_1|^2$ wich is generically different from $J_{st}$.

\vspace{0,5cm}
 
We have  proved in Lemma \ref{lemmod} that when the origin is a point of D'Angelo type  four, then  
$C_{1,\nu}=O\left(|z_1|^{3}\right)+O\left(|z_2|\right)$
and so $(\Lambda_\nu)_*J_\nu=\left(\Lambda_\nu \circ \Psi_{\nu}\right)_*J^\nu$ converges  
to $J_{st}$ when $\nu$ tends to $+\infty$, with: 
\begin{equation*}
\left\{
\begin{array}{lll}  
R_{1,\nu}& = & S_{1,\nu}=0,\\
& & \\
R_{2,\nu}& = & -\alpha_\nu x_1^2+ \alpha_\nu y_1^2,\\
& & \\
S_{2,\nu}& = & -2 \alpha_\nu x_1y_1.\\
\end{array}
\right.
\end{equation*}
%\begin{eqnarray*} 
%R_{2,\nu}& = & S_{2,\nu}=0,\\
%R_{2,\nu}& = & -\alpha_\nu x_1^2+ \alpha_\nu y_1^2, \mbox{ and, }\\
%S_{2,\nu}& = & -2 \alpha_\nu x_1y_1.\\
%\end{eqnarray*} 
 
\vspace{0,3cm}

In case the D'Angelo type of the origin is greater than four, 
we cannot guarantee the convergence of $\tau_\nu \delta_\nu^{-1} C_1^\nu\left(\tau_\nu z_1,\delta_\nu z_2\right)$
when we only remove harmonic terms. So we need to find a more general sequence of  diffeomorphisms $\Psi_\nu$ defined by 
(\ref{eqdif}), (\ref{eqdif1})  and (\ref{eqdif2}) and such that $C_{1,\nu}=O\left(|z_1|^{2m-1}\right)+O\left(|z_2|\right)$.
\vspace{0.7cm}\\
\textbf {Claim.} There are no polynomial $R_{1,\nu}, S_{1,\nu}, R_{2,\nu}$ and $S_{2,\nu}$ such that 
$C_{1,\nu}$ does not contain any order three terms in $x_1$ and $y_1$.
\vspace{0.7cm}

A direct computation leads to:
\begin{eqnarray*}
\alpha_\nu^{-1}\left(J_\nu\right)^3_1\left(z\right)& =  & 
\left(a_2^\nu-a_1^\nu\right)\left(\Psi_\nu^{-1}\left(z\right)\right)x_1-
\left(c_1^\nu+b_2^\nu\right)\left(\Psi_\nu^{-1}\left(z\right)\right)y_1
-y_1\frac{\partial  R_{1,\nu}}{\partial x_1} \\
&&\\
& &-x_1\frac{\partial  R_{1,\nu}}{\partial y_1}
 -x_1\frac{\partial  S_{1,\nu}}{\partial x_1}+y_1\frac{\partial  S_{1,\nu}}{\partial y_1}+
x_1\frac{\partial  R_{1,\nu}}{\partial x_1}\frac{\partial S_{1,\nu}}{\partial x_1}+
y_1\frac{\partial  R_{1,\nu}}{\partial x_1}\frac{\partial S_{1,\nu}}{\partial y_1}\\
&&\\
& &-y_1\frac{\partial  R_{1,\nu}}{\partial y_1}\frac{\partial S_{1,\nu}}{\partial x_1}+ 
x_1\frac{\partial  R_{1,\nu}}{\partial y_1}\frac{\partial S_{1,\nu}}{\partial y_1}
-y_1\left(\frac{\partial S_{1,\nu}}{\partial x_1}\right)^2-y_1\left(\frac{\partial S_{1,\nu}}{\partial y_1}\right)^2\\
&&\\
& &-x_1\frac{\partial  R_{1,\nu}}{\partial x_1}\frac{\partial R_{2,\nu}}{\partial y_2}+
x_1\frac{\partial  S_{1,\nu}}{\partial y_1}\frac{\partial R_{2,\nu}}{\partial y_2}+ 
y_1\frac{\partial  R_{1,\nu}}{\partial y_1}\frac{\partial R_{2,\nu}}{\partial y_2}+\\
&&\\
& & y_1\frac{\partial  S_{1,\nu}}{\partial x_1}\frac{\partial R_{2,\nu}}{\partial y_2}+O\left(|z_1|^4+|z_2|\|z\|\right)
\end{eqnarray*}
and to 
\begin{eqnarray*}
\alpha_\nu^{-1}\left(J_\nu\right)^3_2\left(z\right)& =  & 
\left(b_1^\nu-b_2^\nu\right)\left(\Psi_\nu^{-1}\left(z\right)\right)x_1+
\left(a_1^\nu+a_2^\nu\right)\left(\Psi_\nu^{-1}\left(z\right)\right)y_1+
x_1\frac{\partial  R_{1,\nu}}{\partial x_1}\\
&&\\
& & -y_1\frac{\partial  R_{1,\nu}}{\partial y_1} 
-y_1\frac{\partial  S_{1,\nu}}{\partial x_1}-x_1\frac{\partial  S_{1,\nu}}{\partial y_1}
-x_1\left(\frac{\partial R_{1,\nu}}{\partial x_1}\right)^2-x_1\left(\frac{\partial R_{1,\nu}}{\partial y_1}\right)^2+\\
&&\\
& & y_1\frac{\partial  R_{1,\nu}}{\partial x_1}\frac{\partial S_{1,\nu}}{\partial x_1}+ 
x_1\frac{\partial  R_{1,\nu}}{\partial x_1}\frac{\partial S_{1,\nu}}{\partial y_1}
-x_1\frac{\partial  R_{1,\nu}}{\partial y_1}\frac{\partial S_{1,\nu}}{\partial x_1}+
 y_1\frac{\partial R_{1,\nu}}{\partial y_1}\frac{\partial S_{1,\nu}}{\partial y_1}\\
&&\\
& & -x_1\frac{\partial R_{1,\nu}}{\partial y_1}\frac{\partial R_{2,\nu}}{\partial y_2}-
x_1\frac{\partial S_{1,\nu}}{\partial x_1}\frac{\partial R_{2,\nu}}{\partial y_2}
-y_1\frac{\partial  R_{1,\nu}}{\partial x_1}\frac{\partial R_{2,\nu}}{\partial y_2}+
y_1\frac{\partial  S_{1,\nu}}{\partial y_1}\frac{\partial R_{2,\nu}}{\partial y_2}+\\
&&\\
&& O\left(|z_1|^4+|z_2|\|z\|\right).
\end{eqnarray*}
The only order two terms in $x_1$ and $y_1$ of $\alpha_\nu^{-1}\left(J^\nu\right)^3_1\left(z\right)$ and of 
$\alpha_\nu^{-1}\left(J^\nu\right)^3_2\left(z\right)$ are 
those contained, respectively, in   
$$  
\displaystyle  -y_1\frac{\partial  R_{1,\nu}}{\partial x_1}
-x_1\frac{\partial  R_{1,\nu}}{\partial y_1} -x_1\frac{\partial  S_{1,\nu}}{\partial x_1}+
y_1\frac{\partial  S_{1,\nu}}{\partial y_1}$$  
and 
$$\displaystyle x_1\frac{\partial  R_{1,\nu}}{\partial x_1}-y_1\frac{\partial  R_{1,\nu}}{\partial y_1}
-y_1\frac{\partial  S_{1,\nu}}{\partial x_1}-x_1\frac{\partial  S_{1,\nu}}{\partial y_1} .$$
%\begin{eqnarray*}
% -y_1\frac{\partial  R_{1,\nu}}{\partial x_1}
%-x_1\frac{\partial  R_{1,\nu}}{\partial y_1} -x_1\frac{\partial  S_{1,\nu}}{\partial x_1}+
%y_1\frac{\partial  S_{1,\nu}}{\partial y_1}  & \mbox{, and, }\\  
%x_1\frac{\partial  R_{1,\nu}}{\partial x_1}-y_1\frac{\partial  R_{1,\nu}}{\partial y_1}
%-y_1\frac{\partial  S_{1,\nu}}{\partial x_1}-x_1\frac{\partial  S_{1,\nu}}{\partial y_1}. &
%\end{eqnarray*}
%We need to find $R_{1,\nu}$ and $S_{1,\nu}$ such that there are no such order two terms. 
%But in that case, the limit domain of $\left(\Lambda_\nu\circ \Psi_\nu \Phi^\nu\right) \left(D\cap U\right)$ with respect to 
%the limit structure of $\left(\Lambda_{\nu}\circ \Psi^\nu \circ \Phi_\nu\right)_*J$ has no reason to be Brody hyperbolic.
\vspace{0,4cm}\\
Vanishing these order two terms leads to: 
\begin{equation*}
\left\{
\begin{array}{lll}  
R_{1,\nu}& = & r_{5,\nu}x_1^2-2s_{5,\nu}x_1y_1-r_{5,\nu}y_1^2+r_{1,\nu}x_1^3+r_{2,\nu}x_1^2y_1+r_{3,\nu}x_1y_1^2+r_{4,\nu}y_1^3+\\
& & \\
& & O\left(|z_1|^4+y_2^2+|y_2|\|z\|\right)\\
& & \\
S_{1,\nu}& = & s_{5,\nu}x_1^2+2s_{5,\nu}x_1y_1-s_{5,\nu}y_1^2+s_{1,\nu}x_1^3+s_{2,\nu}x_1^2y_1+s_{3,\nu}x_1y_1^2+s_{4,\nu}y_1^3+\\
& & \\
& & O\left(|z_1|^3+y_2^2+|y_2|\|z\|\right).\\
\end{array}
\right.
\end{equation*}
%\begin{eqnarray*} 
%R_{1,\nu}& = & r_{5,\nu}x_1^2-2s_{5,\nu}x_1y_1-r_{5,\nu}y_1^2+O\left(|z_1|^3+y_2^2+|z_1||y_2|\right), \mbox{ and, }\\
%S_{1,\nu}& = & s_{5,\nu}x_1^2+2s_{5,\nu}x_1y_1-s_{5,\nu}y_1^2+O\left(|z_1|^3+y_2^2+|z_1||y_2|\right).\\
%\end{eqnarray*} 
Then it follows that: 
\begin{eqnarray*}
\alpha_\nu^{-1}\left(J_\nu\right)^3_1\left(z\right)& =  & 
\left(a_2^\nu-a_1^\nu\right)\left(\Psi_\nu^{-1}\left(z\right)\right)x_1-
\left(c_1^\nu+b_2^\nu\right)\left(\Psi_\nu^{-1}\left(z\right)\right)y_1
-y_1\frac{\partial  R_{1,\nu}}{\partial x_1}
\\
&&\\
& & -x_1\frac{\partial  R_{1,\nu}}{\partial y_1}
-x_1\frac{\partial  S_{1,\nu}}{\partial x_1}+y_1\frac{\partial  S_{1,\nu}}{\partial y_1}+O\left(|z_1|^4+|z_2|\|z\|\right),
\end{eqnarray*}
and that
\begin{eqnarray*}
\alpha_\nu^{-1}\left(J_\nu\right)^3_2\left(z\right)& =  & 
\left(b_1^\nu-b_2^\nu\right)\left(\Psi_\nu^{-1}\left(z\right)\right)x_1+
\left(a_1^\nu+a_2^\nu\right)\left(\Psi_\nu^{-1}\left(z\right)\right)y_1
+x_1\frac{\partial  R_{1,\nu}}{\partial x_1}
\\
&&\\
& & -y_1\frac{\partial  R_{1,\nu}}{\partial y_1}
-y_1\frac{\partial  S_{1,\nu}}{\partial x_1}-x_1\frac{\partial  S_{1,\nu}}{\partial y_1}+O\left(|z_1|^4+|z_2|\|z\|\right).
\end{eqnarray*}
\vspace{0,4cm}\\
Since $J^\nu$ satisfies (\ref{eqstrO}), we have: 
\\
\begin{equation*}
\left\{
\begin{array}{lll}  
\left(a_2^\nu-a_1^\nu\right)\left(\Psi_\nu^{-1}\left(z\right)\right)x_1-
\left(c_1^\nu+b_2^\nu\right)\left(\Psi_\nu^{-1}\left(z\right)\right)y_1& = 
&H_{3,\nu}\left(x_1,y_1\right)+\\
&&\\
& & O\left(|z_1|^4+|z_2|\|z\|\right)\\
& & \\
\left(b_1^\nu-b_2^\nu\right)\left(\Psi_\nu^{-1}\left(z\right)\right)x_1+
\left(a_1^\nu+a_2^\nu\right)\left(\Psi_\nu^{-1}\left(z\right)\right)y_1 & = 
&H'_{3,\nu}\left(x_1,y_1\right)+\\
&&\\
& & O\left(|z_1|^4+|z_2|\|z\|\right),\\
\end{array}
\right.
\end{equation*}
\\
%\begin{eqnarray*}
%\left(a_2^\nu-a_1^\nu\right)\left(\Psi_\nu^{-1}\left(z\right)\right)x_1-\left(c_1^\nu+b_2^\nu\right)\left(\Psi_\nu^{-1}\left(z\right)\right)y_1& = &H_{3,\nu}\left(x_1,y_1\right)+O\left(|z_1|^4+|z_2|\|z\|\right)\\
%\left(b_1^\nu-b_2^\nu\right)\left(\Psi_\nu^{-1}\left(z\right)\right)x_1+\left(a_1^\nu+a_2^\nu\right)\left(\Psi_\nu^{-1}\left(z\right)\right)y_1 & = &H'_{3,\nu}\left(x_1,y_1\right)+O\left(|z_1|^4+|z_2|\|z\|\right),\\
%\end{eqnarray*}
where $H_{3,\nu}\left(x_1,y_1\right)$ and $H'_{3,\nu}\left(x_1,y_1\right)$ are real homogeneous polynomials of degree three 
in $x_1$ and $y_1$ which are generically non identically zero. Since we cannot insure the convergence of 
$$\alpha_\nu\tau_\nu\delta_\nu^{-1} H_{3,\nu}\left(\tau_\nu x_1,\tau_\nu y_1\right) =
\alpha_\nu\tau_\nu^4\delta_\nu^{-1} H_{3,\nu}\left(x_1,y_1\right)$$
and 
$$\alpha_\nu\tau_\nu\delta_\nu^{-1} H_{3,\nu}'\left(\tau_\nu x_1,\tau_\nu y_1\right) = 
\alpha_\nu\tau_\nu^4\delta_\nu^{-1} H_{3,\nu}'\left(x_1,y_1\right),$$
we want to cancel polynomials  $H_{3,\nu}\left(x_1,y_1\right)$ and  $H_{3,\nu}'\left(x_1,y_1\right)$ by order three terms in 
$x_1$ and $y_1$ contained in 
$$ \displaystyle -y_1\frac{\partial  R_{1,\nu}}{\partial x_1}-x_1\frac{\partial  R_{1,\nu}}{\partial y_1}
-x_1\frac{\partial  S_{1,\nu}}{\partial x_1}+y_1\frac{\partial  S_{1,\nu}}{\partial y_1}$$
and
$$\displaystyle x_1\frac{\partial  R_{1,\nu}}{\partial x_1}-y_1\frac{\partial  R_{1,\nu}}{\partial y_1}
-y_1\frac{\partial  S_{1,\nu}}{\partial x_1}-x_1\frac{\partial  S_{1,\nu}}{\partial y_1}.$$

Finally, vanishing order three terms in $x_1$ and $y_1$ of 
$\alpha_\nu^{-1}\left(J^\nu\right)^3_1\left(z\right)$ and of $\alpha_\nu^{-1}\left(J^\nu\right)^3_2\left(z\right)$ involve 
the following  system of linear equations: 
\begin{equation*}
\left(\begin{array}{cccccccc} 
3 & 0 & 2 & 0 & 0 & 1 & 0 & 0 \\
3 & 0 & 0 & 0 & 0 & -1 & 0 & 0 \\
0 & 1 & 0 & 0 & 3 & 0 & 0 & 0\\
0 & 2 & 0 & 3 & 0 & 0 & -1 & 0 \\ 
0 & 1 & 0 & 0 & -3 & 0 & -2 & 0 \\
0 & 0 & 1 & 0 & 0 & 0 & 0 & -3 \\
0 & 0 & 1 & 0 & 0 & 2 & 0 & 3 \\
0 & 0 & 0 & 3  & 0 & 0 & 1 & 0 \\ 
\end{array}\right)\left(\begin{array}{ccccc} 
 r_{1,\nu} \\
 r_{2,\nu}\\
 r_{3,\nu}\\
 r_{4,\nu}\\ 
 s_{1,\nu} \\
 s_{2,\nu}\\
 s_{3,\nu}\\
 s_{4,\nu}\\ 
\end{array}\right)
=Y
\end{equation*}
Since this $8\times 8$ system of linear equations is not a Cramer system, it follows that there does not exist, generically, 
polynomials $R_{1,\nu}$ and  $S_{1,\nu}$ such that there are no order 
three term in $x_1$ and $y_1$ in $\left(J^\nu\right)^3_1\left(z\right)$ and $\left(J^\nu\right)^3_2\left(z\right)$. 
\end{proof}

\end{document}